\documentclass[review,10pt]{elsarticle}
\usepackage[left=3cm,top=3cm,right=3cm,bottom=3cm]{geometry}
\usepackage{amsmath, amsthm, amssymb}
\usepackage{lineno,hyperref}
\usepackage{pst-all}
\usepackage{titlesec}
\usepackage{multido}
\usepackage{graphicx}
\usepackage{color}
\usepackage{lipsum}
\usepackage{mathtools}
\usepackage[vcentermath]{youngtab}
\usepackage{young}
\usepackage[all,ps,dvips,graph]{xy}
\usepackage[misc,geometry]{ifsym}

\newcommand*{\Lbullet}{\raisebox{-0.25ex}{\scalebox{.7}{$\bullet$}}}

\let\OLDthebibliography\thebibliography
\renewcommand\thebibliography[1]{
  \OLDthebibliography{#1}
  \setlength{\parskip}{0pt}
  \setlength{\itemsep}{0pt plus 0.3ex}
}

\titleformat{\section} {\normalfont\scshape \large \centering}{ \thesection}{1em}{}
\definecolor{morado}{rgb}{0.5,0,0.5}

\newcommand{\T}{\mathfrak{t}}
\newcommand{\Yla}{Y_n(\lambda)}
\newcommand{\ide}{e(\bs{i}^\lambda)}

\newcommand{\bs}[1]{ \boldsymbol{#1}}
\newcommand{\bip}{ \operatorname{Bip}_1(n)}
\newcommand{\ydl}{Y_n(\lambda)}

\newcommand{\ese}{ {\color{red} s}}
\newcommand{\te}{ {\color{blue} t}}
\newcommand{\al}{ A_{w}}
\newcommand{\B}{ \mathcal{B}}
\newcommand{\X}{ \mathcal{X}}
\newcommand{\C}{ \mathbb{C}}
\newcommand{\blob}{ b_n }
\newcommand{\laba}[1]{ \underline{#1}}
\modulolinenumbers[5]
\newcommand\bi{{\boldsymbol i}}

\input{dibublob.tex}

\newtheorem{teo}{Theorem}[section]
\newtheorem{lem}[teo]{Lemma}

\newtheorem{pro}[teo]{Proposition}
\newtheorem{defi}[teo]{Definition}
\newtheorem{cor}[teo]{Corollary}
\newtheorem{exa}[teo]{Example}
\newtheorem{rem}[teo]{Remark}
\newtheorem{con}[teo]{Conjecture}
\newtheorem{cla}[teo]{Claim}
\newtheorem{condi}[teo]{Conditions}
\newtheorem{algo}[teo]{Algorithm}

\newenvironment{dem}{\noindent \textit{Proof:} }{\quad \hfill $\square$}
\newenvironment{demespecial}{\noindent \textit{\bf{Proof of Theorem} \ref{teo homo sobre}:} }{\quad \hfill $\square$}
\newenvironment{demespecialB}{\noindent \textbf{Proof of Lemma \ref{lema central reflection equal to product of Y's}}. }{\quad \hfill $\square$}

\numberwithin{equation}{section}

\newgray{plomoclaro}{.90}
\newgray{plomooscuro}{.70}
\newgray{blanco}{1}

\begin{document}

\begin{frontmatter}

\title{Blob algebra and two-color Soergel calculus}

\author{Jorge Espinoza\footnote{Instituto de Matem\'aticas y F\'isica, Avenida Lircay s/n, Talca Chile.\newline \hspace*{5mm} email: joespinoza@utalca.cl.}}
\author{David Plaza\footnote{Instituto de Matem\'aticas y F\'isica, Avenida Lircay s/n, Talca Chile.\newline \hspace*{4mm} email: dplaza@inst-mat.utalca.cl.}}

\begin{abstract}
In 2003, Martin and Woodcock noticed a connection between the representation theory of the blob algebra and the Kazhdan--Lusztig polynomials associated with the infinite dihedral group. However, no conceptual explanation for this coincidence has yet been provided. In this study, a possible explanation of this phenomenon is suggested by enunciating a conjecture that relates the endomorphism algebra of Bott--Samelson bimodules to certain subalgebras of the blob algebra obtained by idempotent truncation. Evidence supporting this conjecture is provided.
\end{abstract}


\end{frontmatter}

\section{Introduction}
The purpose of this study is to investigate the relationship  between the blob algebra and the two-color Soergel calculus. The blob algebra $b_n(q,m)$ was introduced by Martin and Saleur in \cite{MSblob} as a two-parameter generalization of the well-known Temperley--Lieb algebra \cite{tl}, whereas the two-color Soergel calculus is a diagrammatic presentation for the category of Soergel bimodules for dihedral groups \cite{elias}. This category controls the behavior of the Kazhdan--Lusztig polynomials.\par

In 2003, Martin and Woodcock noticed a connection between the representation theory of the blob algebra and the Kazhdan--Lusztig polynomials associated with the infinite dihedral group $W=\langle s,t \mbox{ : } s^2=t^2=1 \rangle$. To properly introduce this connection, certain facts regarding the representation theory of $b_n(q,m)$ should be recalled. In what follows, we will work over the field of complex numbers $\C$. We are interested in the  non-semisimple representation theory of the blob algebra. For this reason, and in view of the semisimplicity criterion for $b_n(q,m)$ given in \cite[Section 2]{MWblob00}, a primitive $l$-th root of unity $q\in \C$ and an integer $m$ satisfying $1\leq m <l$ are fixed. For these choices, let $b_n:=b_n(q,m)$.\par
An important feature of $\blob $ is that it is a cellular algebra in the sense of Graham and Lehrer \cite{GLcellularalgebras}. This allows defining a family of $\blob $-modules that are called \emph{cell modules}. Cell modules are parameterized by
$\Lambda_n := \{ \lambda \in \mathbb{Z} \mbox{ :} -n\leq \lambda \leq n\mbox{, } \lambda \equiv n  \mod 2\}$.
Given $\lambda \in \Lambda_n$, the corresponding module is denoted by $\Delta_n(\lambda)$. By the general theory of cellular algebras, each $\Delta_n(\lambda)$ is equipped with a symmetric bilinear form whose radical is a $\blob $-submodule. This allows taking the quotient of $\Delta_n(\lambda)$ and the radical of its bilinear form to obtain a $\blob $-module, which is denoted by $L_n(\lambda)$. The modules $L_n(\lambda)$ are simple for all $\lambda \in \Lambda_n$. Furthermore, each simple $\blob $-module arises in this manner.\par

To calculate the dimensions of the simple $\blob$-modules $L_n(\lambda)$, it suffices to compute the decomposition numbers of $\blob $. Given $\lambda, \mu \in \Lambda_n $,  the \emph{decomposition number} $d_{\lambda , \mu } \in \mathbb{Z}_{\geq 0}$ is defined to be the number of times that $L_n(\mu )$ occurs  as a factor in a composition series of $\Delta_n(\lambda)$. These numbers were calculated by Martin and Woodcock in \cite{MWblob00}. They also observed in \cite[Section 1.1]{MWblob}  that the decomposition numbers of $b_n$ are given by evaluations at $1$ of certain Kazhdan--Lusztig polynomials $h_{x,y}(v) \in \mathbb{Z}[v]$, $x,y \in W$. This observation can be geometrically expressed  as follows. On the real line, vertical lines are drawn at each integer $j$ satisfying $j \equiv -m \mod l$, as is shown in (\ref{dibujo intro}).
\begin{equation}  \label{dibujo intro}
\psscalebox{.55}{\intro }
\end{equation}
The vertical lines in (\ref{dibujo intro}) will be referred to as \emph{walls}, and the connected components of $\mathbb{R} $ with the walls removed will be referred to as \emph{alcoves}.  The alcove containing $0$ is denoted by $\mathcal{A}^0$. $s$ and $t$ are identified with the reflections at $-m$ and $-m+l$, respectively. This defines an action of $W $ on $\mathbb{R}$ and on the set of alcoves. For $x\in W$, let $\mathcal{A}^x:= x \cdot \mathcal{A}^0 $.  It is now assumed that $\lambda, \mu \in \Lambda_n $,   $\lambda \in \mathcal{A}^x$, and $\mu \in\mathcal{A}^y$ for some $x,y \in W$. Then, Martin and Woodcock's observation can be stated as
\begin{equation}  \label{eq intro}
d_{\lambda , \mu} = \left\lbrace
\begin{array}{rl}
h_{x,y}(1),  & \mbox{if }\mu \mbox{ and } \lambda \mbox{ belong to the same orbit;} \\
0,  & \mbox{otherwise.}
\end{array} \right.
\end{equation}
(\ref{eq intro}) gives rise to two natural questions:
\begin{description}
  \item[Q1] Why do polynomials appear in (\ref{eq intro})?
  \item[Q2] Is there a connection between the blob algebra and some object related to the Kazhdan--Lusztig polynomials that explains the combinatorial coincidence in
  (\ref{eq intro})?
\end{description}
An answer to {\textbf{Q1}} was given in \cite{plaza} and \cite{PR-Hblob}. In \cite{PR-Hblob}, it was proven that $\blob $ is endowed with a non-trivial $\mathbb{Z}$-grading by showing that $\blob $ admits a KLR-type presentation, that is, with homogeneous generators $e(\bs{i})$, $y_r$, and $ \psi_s$ (see Theorem \ref{thm isomorphism blob KRL} below). Furthermore, a graded cellular basis for  $\blob $ was constructed, making it a graded cellular algebra in the sense of Hu and Mathas \cite{HMgradedcellular}. The existence of a graded cellular basis allows defining gradings on cell modules and on  simple modules. In this graded setting, the decomposition numbers become Laurent polynomials. They are called \emph{graded decomposition numbers} and are denoted by $d_{\lambda , \mu }(v) \in \mathbb{Z}[v,v^{-1}]$. In \cite{plaza}, the graded decomposition numbers of $\blob $ were calculated explicitly.  On the other hand, by taking advantage of the fact that there exists a closed formula for the Kazhdan--Lusztig polynomials associated with $W$, it can be concluded that
\begin{equation}  \label{eq intro dos}
d_{\lambda , \mu}(v) = \left\lbrace
\begin{array}{rl}
h_{x,y}(v),  & \mbox{if }\mu \mbox{ and } \lambda \mbox{ belong to the same orbit;} \\
0,  & \mbox{otherwise,}
\end{array} \right.
\end{equation}
for all $x,y \in W$ and $\lambda, \mu \in \Lambda_n $ such that $\lambda \in \mathcal{A}^x$ and $\mu \in\mathcal{A}^y$. Of course, (\ref{eq intro}) is obtained from (\ref{eq intro dos}) by evaluating at $v=1$. Thus, (\ref{eq intro dos}) may be considered an improvement of (\ref{eq intro}).\par
It is an usual phenomenon in mathematics that equalities between numbers or polynomials are the tip of an iceberg of a deeper fact. The purpose of this paper is to try to figure out what  such an ``iceberg'' is for the equality (\ref{eq intro dos}). To give an accurate answer to the above, some notation is required.\par
A homogeneous idempotent $\ide \in  \blob$ can be associated with each $\lambda \in \Lambda_n $ (see Definition \ref{defi sequence residues} below). Then, the idempotent truncation $b_n(\lambda ):= \ide b_n \ide $ of $b_n$ is considered. These subalgebras were studied in \cite{plaza} to compute the graded decomposition numbers of $b_n$. There, it was shown that $b_n(\lambda )$ is a graded cellular algebra with graded cellular basis formed by the elements of the graded cellular basis of $b_n$ that belong to $b_n(\lambda)$. Furthermore, as $b_n(\lambda)$ is obtained as an idempotent truncation of $b_n$, the graded decomposition numbers of $b_n$ and $b_n(\lambda)$ coincide. On the other hand, given $w\in W$, a $\mathbb{Z}$-graded $\C$-algebra $\al$ is considered. This algebra is defined diagrammatically in Section 2. It corresponds to the diagrammatic version, given in \cite{elias}, of the endomorphism algebra of the  Bott--Samelson bimodule associated with $w$, with the left action trivialized. The algebra $\al$ is also a graded cellular algebra \cite[Theorem 4.1]{plazauno} whose graded cellular basis is the Double Leaves Basis introduced by Libedinsky in \cite{nico}. In this case, the graded decomposition numbers are given by the Kazhdan--Lusztig polynomials associated with $W$ \cite[Theorem 4.8]{plazauno}. Then, we have the following conjecture.
\begin{con} \label{conjecture one}
  Let $\lambda \in \Lambda_n$ and $w\in W$. It is assumed that $\lambda \in \mathcal{A}^w$. Then, $ b_n(\lambda )\cong \al$ as graded cellular algebras.
\end{con}
To be isomorphic as graded cellular algebras implies that an algebra isomorphism between $ b_n(\lambda )$ and $ \al $ exists, and that this isomorphism is compatible with the graded cellular structure. Specifically, the isomorphism maps the graded cellular basis of $ b_n(\lambda )$ to the graded cellular basis of $\al$. In the final section, the expected isomorphism is explicitly described by establishing a correspondence between the graded cellular bases.  Of course, Conjecture \ref{conjecture one} implies (\ref{eq intro dos}).
\par
The blob algebra is the smallest example of a family of  Hecke algebra quotients known as \emph{Generalized blob algebras} introduced in  \cite{MWblob}. In \cite{bowman} a generalization of (\ref{eq intro dos}) was proven for the aforementioned algebras. Concretely, Bowman demonstrated that generalized blob algebras are graded cellular algebras and that (when the ground field is of characteristic zero) their graded decomposition numbers coincide with affine Kazhdan-Lusztig polynomials of type $A$. It should be mentioned that Bowman's graded cellular basis still makes sense  when generalized blob algebras are defined over a field of positive characteristic, but the graded decomposition number are no longer affine Kazhdan-Lusztig polynomials of type $A$. Two natural questions arise in this case:
\begin{description}
  \item[Q3] Is there a conceptual explanation for the generalization of (\ref{eq intro dos}) obtained by Bowman? 
   \item[Q4] What are the graded decomposition numbers for the generalized blob algebras in positive characteristic? 
\end{description}

A conjectural answer for {\bf{Q3}} and {\bf{Q4}} is given in \cite{nicoplaza}, where Libedinsky and the second author introduced the so-called \emph{Categorical Blob vs Soergel Conjecture}. For $l\geq 2$ let $W_{l}$ be the affine Weyl group of type $\tilde{A}_{l-1}$.  One of the consequences of the above conjecture is the fact that for any element $w\in W_l$ there exists a reduced expression $\underline{w}$ of $w$  such that, the endomorphism ring of the Bott-Samelson bimodule associated to $\underline{w}$ is isomorphic to certain idempotent truncation of a generalized blob algebra of level $l$. 
    So that the Categorical Blob vs Soergel Conjecture implies the natural generalization of Conjecture \ref{conjecture one}  to higher levels. Another consequence of the \emph{Categorical Blob vs Soergel Conjecture} is the fact that in characteristic $p$, the graded decomposition numbers of the generalized blob algebras are given by $p$-Kazhdan-Lusztig polynomials \cite{jenwil}. Summing up, the possibility that open questions in the theory of Soergel bimodules may be accessible by studying (generalized) blob algebras makes Conjecture \ref{conjecture one} (and its generalizations) interesting objects of study.\par
 The main result of this paper is to provide evidence supporting Conjecture \ref{conjecture one}. The main obstacle in order to prove such a conjecture is the lack of a presentation  for the two algebras involved in it. In this context in order to provide evidence supporting Conjecture  \ref{conjecture one} it is natural to try to find isomorphic subalgebras of the algebras considered in Conjecture \ref{conjecture one} that, on the one hand, are large enough so that an isomorphism between them gives real evidence in favor of the conjecture, and that on the other hand are  small enough so that a presentation for them exists. In this setting we consider $\Yla \subset b_n(\lambda )$ and $\mathcal{DL}_k \subset \al$, where $k=l(w)$ and $\lambda \in \mathcal{A}^w$. The above subalgebras are introduced in Definition \ref{defi dot line diagram} and Definition \ref{defi Yla}. We will prove that these subalgebras are isomorphic. To this end, a third algebra $DL_k$, which is called the \emph{Dot--Line algebra}, is first defined. Then, it is separately proved that both $\Yla$ and $\mathcal{DL}_k$  are isomorphic to $DL_k$. We stress that $\dim DL_k =2^k$, so that the algebras  $\Yla $ and $\mathcal{DL}_k $ are big enough so that the isomorphism between them really gives  evidence in favor of Conjecture (\ref{conjecture one}). On the other hand, $\Yla$ is by definition the algebra generated  by $\{ \ide, y_1\ide, \ldots , y_n\ide \}$. The KLR-generators $\{y_i\}$ are nilpotent versions of the classical  Jucys-Murphy elements. In this setting,  $\Yla$ can be thought  as the Gelfand-Tsetlin  subalgebra of $b_n(\lambda)$. Which makes $\Yla$ an interesting object of study by its own.\par

The strategy for proving the aforementioned isomorphisms is quite transparent. Namely, $DL_k$ is defined as an algebra with generators and relations, and then these relations are verified in $\Yla$ and $\mathcal{DL}_k$. Finally, the isomorphism is established by a dimension argument. It is easy to prove the isomorphism between $DL_k$ and $\mathcal{DL}_k$. In this case, to verify that the relations defining $DL_k$ are satisfied in 
$\mathcal{DL}_k$ is a straightforward exercise, as well as to calculate $\dim_\mathbb{C} \mathcal{DL}_k$. By contrast, to verify the isomorphism between $DL_k$ and $\Yla$ is not an easy task, and indeed a major part of the study is dedicated to this. The first obstacle in obtaining the desired isomorphism is that it is unclear as to which elements should be the images of the generators of $DL_k$. Even after the correct elements in $\Yla$ have been found, the fact that they satisfy the relations defining $DL_k$ does not follow immediately. In this setting, the homogeneous presentation of $b_n(q,m)$ (see Theorem \ref{thm isomorphism blob KRL}) should be considered, which (as the reader can convince himself/herself by looking at such a presentation) is highly complicated. Instead of using this presentation, a variant of the diagrammatic calculus introduced by Khovanov and Lauda in \cite{KLdiagram} is utilized. The diagrammatic approach allows performing cumbersome algebraic computations more transparently. However, it should be noted that this approach still has certain subtleties, and to achieve the desired goal, delicate arguments involving these diagrams should be utilized, which we believe are of independent interest. Finally, to calculate the dimension of $\Yla$ is also difficult. In fact, it does not directly follow from the definition of $\Yla$ even that $\Yla \neq 0$. This obstacle is overcome by resorting to the graded cellular basis of $b_n$ constructed in \cite{PR-Hblob}.
\par

The structure of the paper is as follows. In the next section, the two-color Soergel calculus is introduced. Specifically, with each $w\in W$, a $\mathbb{Z}$-graded $\mathbb{C}$-algebra $\al $ is associated that is defined diagrammatically. The construction of the Double Leaves Basis for $\al $ is also reviewed. In section $3$, the Dot--Line algebra is defined and is realized as a subalgebra of $\al $.  In section $4$, the blob algebra with its graded presentation is introduced, as well as the diagrammatic setting for working with it. Section $5$ provides the combinatorial background to perform calculations with Khovanov--Lauda diagrams. In section $6$, the behavior of the KLR generators $y_1, \ldots ,y_n \in b_n(\lambda)$ is studied.  In section $7$, the isomorphism between the Dot--Line algebra and $\Yla$ is established. Finally, in the last section, a graded vector space isomorphism between $\al $ and $b_n(\lambda)$ is provided. Although this isomorphism is  combinatorially constructed, we firmly believe that it is the algebra isomorphism that would allow proving Conjecture \ref{conjecture one}.

\section*{Acknowledgments}
We would like to thank to the anonymous referee for her/his comments and suggestions that helped us  improve the text. The first author was partially  supported by Beca Doctorado Nacional 2013-CONICYT, 21130109. The second author was partially supported by PAI-CONICYT-Concurso Nacional de Inserci\'on en la Academia 2015, 79150016, and Fondecyt Iniciaci\'on project 11160154.

\section{The two-color Soergel calculus}
In this section, the two-color Soergel calculus associated with the infinite dihedral group is introduced by following the diagrammatic approach given by Elias  in \cite{elias}. For the sake of brevity, Soergel bimodules will not be mentioned. For an explanation about the relationship between Soergel bimodules and the diagrams that we are going to  define, the reader is referred to  \cite{EliasKhovanov}, \cite{elias}, or \cite{ewsoergel}. Furthermore, for simplicity, features of the Soergel calculus that do not arise in the two-color setting are ignored. Let $(W ,S)$ be the infinite dihedral group. That is, $S=\{  {\color{red} s}, {\color{blue} t} \}$ and $W$ has a presentation
\begin{equation} \label{presentation W}
W= \langle  {\color{red} s}, {\color{blue} t} \mbox{ ; }   {\color{red} s}^2 ={\color{blue} t}^2=e \rangle .
\end{equation}
The elements in $S $ are called \emph{simple reflections} or \emph{colors}. Given a positive integer $k$, let
\begin{equation}
k_{ {\color{red} s}}:=  \underbrace{{\color{red} s} {\color{blue} t}  {\color{red} s} \ldots }_{k\mbox{-times}}  \qquad k_{ {\color{blue} t}}:=  \underbrace{ {\color{blue} t} {\color{red} s} {\color{blue} t}  \ldots }_{k\mbox{-times}}
 \end{equation}
Moreover, $0_{ {\color{red} s}}=0_{ {\color{blue} t}} =e$. It is clear from (\ref{presentation W}) that each element in $W$ is of the form $k_{ {\color{red} s}}$ or $ k_{ {\color{blue} t}}$ for some integer $k\geq 0$. Soergel calculus depends heavily on a \emph{realization}  $\mathfrak{h}$ of $(W,S)$, a concept that was defined by Elias and Williamson in \cite[Section 3.1]{ewsoergel}. Heuristically, a realization of $(W ,S)$ is a representation of $W$ that satisfies certain technical conditions. In this study, the \emph{geometric representation} of $W$ defined over the complex numbers will be considered \cite[Section 5.3]{hum}. Let $\mathfrak{h}:= \mathbb{C}\alpha^\vee_\ese \oplus \mathbb{C}\alpha^\vee_\te $. The assignment

\begin{equation}  \label{matrix realization}
\ese \rightarrow
\left(
\begin{array}{cc}
-1 & 2 \\
0  & 1
\end{array}
\right)
\qquad
\te \rightarrow
\left(
\begin{array}{cc}
1 & 0 \\
2 & -1
\end{array}
\right)
\end{equation}
defines a $\mathbb{C}$-representation of $W$ over $\mathfrak{h}$. Certain elements $\alpha_\ese , \alpha_\te \in \mathfrak{h}^\ast$ are now fixed. They are determined by the following rules:

\begin{equation}  \label{equations realization}
\begin{array}{lll}
\alpha_{\ese } ( \alpha^\vee_\ese )= 2,   & \mbox{ } & \alpha_\te (\alpha^\vee_\ese )= -2 ;    \\
\alpha_\ese (\alpha^\vee_\te ) = -2 , &    & \alpha_\te (\alpha^\vee_\te )= 2  .
\end{array}
\end{equation}
Let $R := S(\mathfrak{h}^\ast) = \oplus_{i\geq 0} S^i(\mathfrak{h}^\ast)$ be the symmetric algebra of $\mathfrak{h}^\ast$. This is a graded $\mathbb{C}$-algebra with the usual grading doubled, i.e., $\deg (\mathfrak{h}^\ast)=2$. $W $ acts on  $\mathfrak{h} $. Hence, it also acts on $\mathfrak{h}^\ast$ via the contragradient representation, and this action extends to $R$. The \emph{Demazure operators} $\partial_\ese , \partial_\te :R \rightarrow R(-2)$ are defined by

\begin{equation}
\partial_\ese (f) = \frac{f-\ese f}{\alpha_\ese},  \qquad \qquad \partial_\te (f) = \frac{f-\te f}{\alpha_\te}.
\end{equation}

\begin{defi}
A two-color Soergel graph (or simply an $S$-graph for short) is a finite and  decorated graph with its boundary properly embedded into the planar strip $\mathbb{R}\times [0,1]$. The edges in an $S$-graph are colored by $\ese $ and $\te $. The vertices in this graph are of two types (see (\ref{vertices})):
\begin{itemize}
 \item  Univalent vertices (dots). These have degree $+1$.
 \item Trivalent vertices, where all three incident edges have the same color. These have degree $-1$.
 \end{itemize}
 \begin{equation}  \label{vertices}
 \psscalebox{.7}{\vertices }
\end{equation}
Furthermore, an $S$-graph may have its regions (the connected components of the complement of the graph in the strip $\mathbb{R}\times [0,1] $) decorated by boxes labeled by homogeneous $f\in R$. The \emph{degree} of an $S$-graph is defined as the sum of the degrees of each vertex plus the sum of the degrees of the polynomials inside each box.
\end{defi}

Here is an example of an $S$-graph:
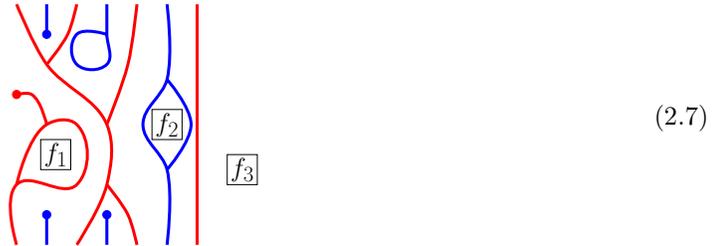
\begin{equation}  \label{ejemploSoergel}
 \psscalebox{.4}{\ejemploSoergel }
\end{equation}
where $f_i \in R$ are homogeneous polynomials. In this example, there are eight trivalent vertices and four dots. Consequently, the degree of the $S$-graph in (\ref{ejemploSoergel}) is $-4+\sum_{i=1}^3 \deg f_i$. \par

The points where an edge touches the boundary of the strip $\mathbb{R}\times [0,1] $  are called \emph{boundary points}. The boundary points of an $S$-graph on $\mathbb{R}\times \{0\}$ and on $\mathbb{R}\times \{1\}$ provide two sequence of colored points. These sequences are called the \emph{bottom boundary} and \emph{top boundary}, respectively. In this study, only $S$-graphs satisfying the following conditions are considered:
\begin{itemize}
\item The bottom boundary and the top boundary coincide.
\item The bottom boundary (and hence the top boundary) is alternating. Therefore, this sequence can be identified with a reduced expression of an element of $W $.
\end{itemize}

For instance, the diagram in (\ref{ejemploSoergel}) satisfies the above conditions. In this case, the bottom/top boundaries are identified with $7_\ese \in W$. A graded $\mathbb{C}$-algebra is now associated with each element $w\in W $. It will be denoted by $A_w$ and will be referred to as the \emph{endomorphism algebra of} $w$.

\begin{defi} \label{defin endo BS}
Given $w\in W$, the endomorphism algebra of $w$ is denoted by $A_w$ and is defined to be the $\mathbb{C}$-vector space generated by isotopy classes of $S$-graphs with bottom boundary and top boundary identified with $w$ modulo the following local relations:

\begin{equation} \label{sgraphA}
\psscalebox{.7}{\relSA }\quad = \quad \psscalebox{.7}{\relSB }
\end{equation}

\begin{equation}   \label{sgraphB}
\psscalebox{.7}{\relSC}\quad = \quad \psscalebox{.7}{\relSD }
\end{equation}

\begin{equation}   \label{sgraphC}
\psscalebox{.7}{\relSF}\quad = \quad \psscalebox{.7}{\relSG }
\end{equation}

\begin{equation}   \label{sgraphD}
\psscalebox{.7}{\relSH}\quad =\quad \psscalebox{.7}{\relSI}\quad + \quad \psscalebox{.7}{\relSJ}
\end{equation}

\begin{equation}     \label{sgraphE}
\psscalebox{.7}{\relSE} \mbox{ } = 0
\end{equation}

\begin{equation}  \label{rel muere a la izquierda}
\psscalebox{.7}{\relSK } \mbox{ }= 0,\qquad \mbox{ if } \deg (f)>0.
\end{equation}

Relations  {\rm(\ref{sgraphA})--(\ref{sgraphE})} hold if red is replaced by blue. Relation (\ref{rel muere a la izquierda}) implies that if a diagram is decorated in its leftmost region by a homogeneous polynomial $f \in R$ of positive degree, then it is equal to zero. Finally, the multiplication on $A_w$ is determined in $S$-graphs by concatenation.
\end{defi}

\begin{rem}\rm
 Let $\underline{w}$ be the unique reduced expression of an element $w\in W $. If  $\mbox{BS}(\underline{w})$ denotes the Bott-Samelson bimodule associated to $\underline{w}$ then  $
		A_w\simeq \mathbb{C} \otimes_{R}\mbox{End} (\mbox{BS}(\underline{w} ))$.
	In other words, $\mbox{End} (\mbox{BS}(\underline{w} ))$ can be recovered from $A_w$ by dropping (\ref{rel muere a la izquierda}).
\end{rem}

This section is concluded by reviewing the construction of the \emph{Double Leaves Basis} for $\al$. This basis was introduced by Libedinsky \cite{nico} at the level of generality of an arbitrary Coxeter system.  Elias and Williamson transferred this basis into the diagrammatic setting \cite{ewsoergel}.  As before, the focus will be only on the two-color Soergel calculus associated with the fixed realization (\ref{matrix realization})--(\ref{equations realization}) of $W$; therefore, features of the Double Leaves Basis that do not arise in this setting will be ignored.
\par
Henceforth, an element $w=s_1\ldots s_n \in W$ is fixed, where $s_i \in \{\ese , \te \}$.  To describe the Double Leaves Basis for $\al $, a perfect binary tree $\mathbb{T}_w$ is first defined with nodes decorated by $S$-graphs. $\mathbb{T}_w$ is constructed by induction on the depth of the nodes. In depth one, $\mathbb{T}_w$ is given by
\begin{equation}
  \psscalebox{.5}{ \arbR }
\end{equation}
Let $1<k \leq n$, and it is assumed that $\mathbb{T}_w$ has already been constructed in depth $j$ for all $1\leq j <k$. Let $N$ be a node of depth $k-1$ decorated by the $S$-graph $f_N$. Then, we have two cases.

\begin{enumerate}
\item If $s_k$ is equal to the rightmost color of the top sequence of $f_N$, then the child nodes of $N$ are decorated as follows:
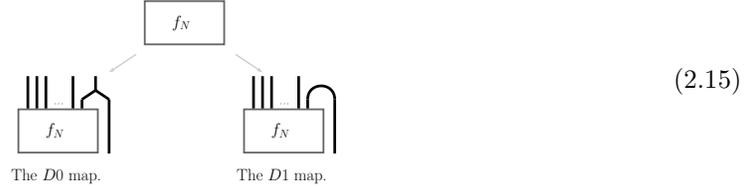
\begin{equation}
  \psscalebox{.5}{ \arbN }
\end{equation}
\item If $s_k$ is different from the rightmost color of the top sequence of $f_N$, or if the top sequence of $f_N$ is the empty sequence, then the child nodes of $N$ are decorated as follows:
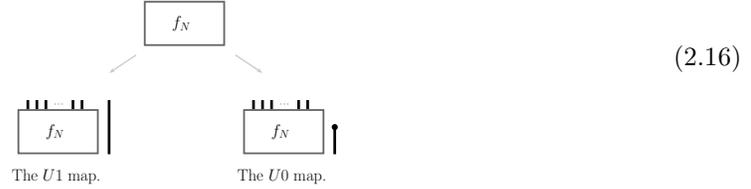
\begin{equation}
  \psscalebox{.5}{ \arbNN }
\end{equation}
\end{enumerate}

\medskip
This completes the construction of $\mathbb{T}_w$. The $S$-graphs that decorate the leaves nodes of $\mathbb{T}_w$ are called \emph{light leaves} of $w$. By construction, it is straightforward to see that the bottom sequence of each light leaf of $w$ is equal to $w$, whereas the top sequence of a light leaf of $w$ can be identified with an element of $W$, that is, it is an alternating sequence of colors. Given $x\in W$, $\mathbb{L}_w(x)$ denotes the set of all light leaves with top sequence equal to $x$. Moreover, $\mathbb{L}_w$ denotes the set of all light leaves of $\mathbb{T}_w$. Given $f\in \mathbb{L}_w$, $f^a$ denotes the $S$-graph obtained from $f$ by applying a vertical flip. Then, let
\begin{equation}
\mathbb{DL}_w = \{  f_1^a \circ f_2 \mbox{ }\vert \mbox{ } f_1,f_2 \in \mathbb{L}_w(x)\mbox{, } x\in W \},
\end{equation}
where $\circ $ denotes concatenation of diagrams. The elements of $\mathbb{DL}_w $ are called \emph{Double Leaves} of $w$. It should be noted that the bottom sequence and the top sequence of each double leaf are equal to $w$. Therefore, $\mathbb{DL}_w$ is a subset of $ \al$. The following theorem is due to Libedinsky \cite{nico}. It was reformulated in the diagrammatic setting by Elias and Williamson \cite{ewsoergel}.

\begin{teo}
Let $w\in W$. Then, $\mathbb{DL}_w$ is a basis of $A_w$.
\end{teo}

\begin{exa}
The construction of the Double Leaves Basis of $A_w$ will be illustrated when $w=\ese \te \ese$. In this case, $\mathbb{T}_w$ is given by
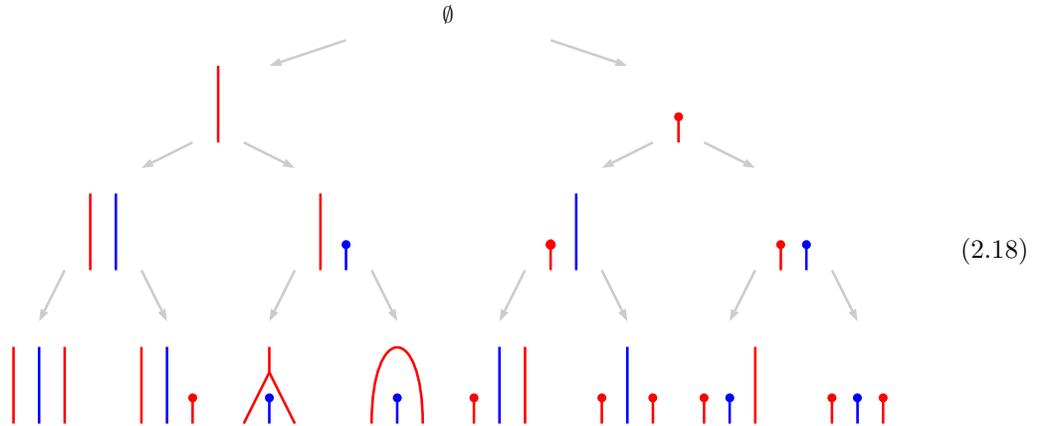
\begin{equation}
\treeA
\end{equation}
 Hence, the Double Leaves Basis of $A_w$ is given by
$$\LLbase$$
\end{exa}

\begin{rem} \label{dot DL is diferent from zero}
Let $w\in W$. It is easy to see that the diagram
\begin{equation}  \label{dibujo double dots in DL is diferent from zero}
\psscalebox{.7}{\lineapunto }
\end{equation}
belongs to the Double Leaves Basis of $\al $. In particular, this diagram is different from zero.
 \end{rem}

\section{The Dot--Line algebra $DL_n$}

In this section, for each positive integer $n$, a commutative algebra, denoted by $DL_n$, is defined and referred to as the \emph{Dot--Line algebra on $n$-strings} (the reason for this terminology will become obvious shortly). $DL_n$ is first defined as an algebra with generators and relations, and subsequently it is shown that $DL_n$ can be embedded as a subalgebra of $A_w$ for all $w\in W$ with $l(w)=n$.

\begin{defi}
The Dot--Line algebra on $n$-strings $DL_n$ is the associative $\mathbb{C}$-algebra generated by $\{1,X_1,\ldots , X_n\} $ subject to the relations
\begin{align}
 X_iX_j&= X_jX_i,  &  & \mbox{ \rm for all } 1\leq i, j \leq n;  \label{dotlineone}  \\
X_i ^2&=  -2X_i\left(\sum_{j=1}^{i-1}X_j \right),  &  &  \mbox{ \rm for all } 2\leq i\leq n; \label{dotlinetwo} \\
 X_1^2 &=0.  & &  \label{dotlinethree}
\end{align}
\end{defi}

Given a sequence $\alpha=(\alpha_1,\ldots , \alpha_n) \in \left( \mathbb{N}_0 \right)^n$, the monomials $X^\alpha \in DL_n$ are defined by  $X^\alpha := X_1^{\alpha_1 } \ldots X_n^{\alpha_n}  $.

\begin{lem} \label{lema span set for Dot--Line}
The algebra $DL_n$ is spanned by $B_n:=\{ X^\alpha \mbox{ } | \mbox{ } \alpha \in \{0,1\}^n  \}$. In particular, $\dim DL_n \leq 2^{n}$.
\end{lem}
\begin{dem}
Let $\beta= (\beta_1 , \ldots , \beta_n) \in \left( \mathbb{N}_0 \right)^n$. It will be shown that $X^{\beta}$ is a linear combination of elements in $B_n$. If $\beta  \in \{0,1\}^n $, there is nothing to prove. Thus, it can be assumed that $\beta_i \geq 2$ for some $1\leq i \leq n$. The proof is by induction on $k_\beta:= \max \{ i \in \{1, \ldots , n\} \mbox{ } | \mbox{  } \beta_i \geq 2 \} $. If $k_\beta=1$, then $X^\beta =0$ by (\ref{dotlinethree}). It is now assumed that $k_\beta > 1$, that is,  $\beta_i \in \{0,1\}$ for $i> k_\beta$. Then, by  repeated applications of (\ref{dotlinetwo}) we have
\begin{align}
  X^{\beta} &= X_1    \ldots   X_{k_\beta-1}^{\beta_{k_\beta-1}}   X_{k_\beta}^{\beta_{k_\beta}} \ldots X_n^{\beta_n}  \\
            &= X_1    \ldots   X_{k_\beta-1}^{\beta_{k_\beta-1}} \left( -2 \sum_{j=1}^{k(\beta)-1} X_j    \right)^{\beta_{k_\beta}-1}    X_{k_\beta} \ldots X_n^{\beta_n}.
\end{align}
Therefore, $X^{\beta}$ is the sum of a linear combination of elements $X^{\alpha}\in B_n$  and a linear combination of elements $X^\gamma \not \in B_n$ with $k_{\gamma} < k_\beta$. Thus, by induction, $X^{\beta}$ is a linear combination of elements in $B_n$.
\end{dem}

\medskip
It will now be shown that $B_n$ is a basis of $DL_n$. An element $w\in W$ is fixed, and its length is denoted by $n$. To prove that $B_n$ is a basis of $DL_n$, it will be useful to realize $DL_n$ as a diagram subalgebra of $A_{w}$. For $1\leq i \leq n$, $\mathcal{X}_i$ denotes the $S$-graph with top boundary and bottom boundary equal to $w$ given by
\begin{equation}   \label{dotlineejemplo}
\psscalebox{.3}{\dlA } \mbox{ , }
\end{equation}
where the dots are located at the $i$-th position. By abuse of notation, $\mathcal{X}_i$ will denote the corresponding  equivalence class in $\al$.  It is easy to see that $\mathcal{X}_i\mathcal{X}_j=\mathcal{X}_j\mathcal{X}_i $ for all $1\leq i,j\leq n$.

\begin{defi}  \label{defi dot line diagram}
Given a sequence $\alpha=(\alpha_1,\ldots , \alpha_n) \in \left( \mathbb{N}_0 \right)^n$, let
\begin{equation}
 \mathcal{X}^\alpha = \prod_{i=1}^n \mathcal{X}_i^{\alpha_i}
\end{equation}
and $\mathcal{B}_n:= \{\mathcal{X}^\alpha \mbox{ } | \mbox{ } \alpha \in \{0,1\}^n \}$. Moreover, $\mathcal{DL}_n$ denotes the subspace of $\al $ spanned by $\mathcal{B}_n $.
\end{defi}

\begin{lem} \label{lem dot line presentation diagrams}
  Let $w\in W$ with $l(w)=n$. Then, $\mathcal{DL}_n$ is a subalgebra of $\al $. Furthermore, there exists a surjective algebra homomorphism $ \eta : DL_n \rightarrow \mathcal{DL}_n$ determined by $\eta (X_i) = \mathcal{X}_i$ for all $1\leq i \leq n$.
\end{lem}

\begin{dem}
 By (\ref{equations realization}), we have $\alpha_\ese + \alpha_\te =0$. Hence, (\ref{sgraphC}) implies
 \begin{equation}\label{relation Double dots}
 \psscalebox{.7}{  \relDoubleDotsA} \quad =- \quad \psscalebox{.7}{  \relDoubleDotsB} \mbox{ .}
 \end{equation}
 Moreover, by combining (\ref{sgraphC}) and (\ref{sgraphD}), it is not difficult to see that
\begin{equation}\label{relation Double dots dos}
 \psscalebox{.7}{ \DDC  } \quad = \quad  \psscalebox{.7}{ \DDD  } \quad -2 \quad \psscalebox{.7}{ \DDE } \mbox{ .}
\end{equation}
Relation (\ref{relation Double dots dos}) has an obvious counterpart that is obtained by switching the colors. By applying (\ref{rel muere a la izquierda}), (\ref{relation Double dots}), and (\ref{relation Double dots dos}), it is straightforward to see that the product of two elements in $\mathcal{B}_n$ is a linear combination of elements in $\mathcal{B}_n$. Therefore, $\mathcal{DL}_n$ is a subalgebra of $\al$.

\par
To prove the existence of the homomorphism $\eta $, it is enough to show that the elements $\mathcal{X}_i$
satisfy the relations (\ref{dotlineone})--(\ref{dotlinethree}). Relation (\ref{dotlineone}) is clear, and (\ref{dotlinethree}) follows from (\ref{rel muere a la izquierda}). Furthermore, by repeated application of (\ref{relation Double dots dos}) and (\ref{rel muere a la izquierda}), (\ref{dotlinetwo}) is obtained, which proves the existence of $\eta$. Finally, the surjectivity  of $\eta$  follows directly from the definitions.
\end{dem}

\medskip
It will be proved that the homomorphism $\eta$ in Lemma \ref{lem dot line presentation diagrams} is in fact an isomorphism. To this end, an order $ \rhd  $ on $\B_n$ should be defined as follows: $\X^{\alpha}\rhd \X^{\beta}$ if and only if $\deg \X^{\alpha} > \deg \X^{\beta}$, or, $\deg \X^{\alpha} =\deg \X^{\beta}$ and $\alpha \succ \beta$, where $\succ$ denotes the usual lexicographical order. 

\par
Given  $\alpha=(\alpha_1, \ldots , \alpha_n) \in \{0,1\}^n$,   $\widehat{\alpha} =(\widehat{\alpha}_1, \ldots , \widehat{\alpha}_n )\in \{0,1\}^n$ is defined as the sequence determined by $\widehat{\alpha}_i  = 1- \alpha_i  $ for all $1\leq i \leq n$. It is clear that $\X^{\alpha}\X^{\widehat{\alpha}}=\X^{(1,1,\ldots,1)}$.

\begin{lem}  \label{lemma mata segun el orden}
Let $w\in W$ and $n=l(w)$. Moreover, let $\alpha=(\alpha_1,\ldots,\alpha_n) \in \{0,1\}^n$. Then $\X^{\beta}\X^{\widehat{\alpha}}=0$ for all $\X^{\beta} \in \B_n$ such that $\X^{\beta}\rhd \X^{\alpha}$.
\end{lem}

\begin{dem}
By the construction of the Double Leaves Basis, the maximal possible degree for an element of $\al $ is $2n$. If $\deg(\X^{\beta})> \deg(\X^{\alpha})$, then $\deg (\X^{\beta}\X^{\widehat{\alpha}}) >2n$. Consequently, $\X^{\beta}\X^{\widehat{\alpha} }=0$, and the lemma is proved in this case. Therefore, it may be assumed that $ \deg (\X^{\beta}) = \deg (\X^{\alpha})$. By the definition of the order $\rhd$, we have $\beta\succ\alpha$. Let $\beta = (\beta_1,\ldots , \beta_n)$. As $\beta\succ\alpha$, there exists an index $i$ such that  $\alpha_i=0$, $\beta_i=1$, and $\alpha_j=\beta_j$ for all $j=1,\ldots,i-1$. Then, the diagram $\X^{\beta}\X^{\widehat{\alpha}}$ has the form
\begin{equation}
\lemA
\end{equation}
The lemma now follows form (\ref{sgraphC}) and (\ref{rel muere a la izquierda}).
\end{dem}

\begin{teo} \label{teo base de Dot line}
  Let $w\in W$ with $n = l(w)$. Then, $\B_n$ is a linearly independent subset of $A_w$. In particular, $\dim \mathcal{DL}_n = 2^n$. Furthermore, $DL_n\cong \mathcal{DL}_n $.
\end{teo}

\begin{dem}
It is enough to prove the linear independence of $\B_n$. The other assertions follow at once from Lemma \ref{lema span set for Dot--Line} and Lemma \ref{lem dot line presentation diagrams}. It is assumed that
\begin{equation}\sum_{\alpha\in\{0,1\}^n}r_{\alpha}\X^{\alpha}=0,\label{eclema2}\end{equation}
for some $r_{\alpha}\in \mathbb{C}$. Then, $r_\alpha =0$, for all $\alpha \in \{0,1 \}^n$.  Indeed, if this is not the case, there exists $\X^{\beta}\in \B_n$ that is minimal with respect to $\rhd$ satisfying $r_{\beta}\neq 0$. By multiplying  (\ref{eclema2}) by $\X^{\widehat{\beta}}$ and then by using  Lemma \ref{lemma mata segun el orden}, we have  $r_{\beta}\X^{(1,1,\ldots,1)}=0$. It is not difficult to see that $\X^{(1,1,\ldots,1)}$ corresponds to the diagram in (\ref{dibujo double dots in DL is diferent from zero}). Therefore,  Remark \ref{dot DL is diferent from zero} implies that $r_{\beta}=0$, which contradicts the choice of $\X^\beta$.
\end{dem}

\section{The blob algebra}
In this section, the blob algebra is defined. This algebra was introduced by Martin and Saleur \cite{MSblob} as a generalization of the Temperley--Lieb algebra. The blob algebra can be realized as a diagram algebra with a diagrammatic basis given by blobbed Temperley--Lieb diagrams, which explains its name.

\par
Although in this study the graded presentation of the blob algebra introduced by the second author and Ryom-Hansen \cite{PR-Hblob} will be considered, for completeness, the classical definition of the blob algebra is also provided (see \cite[Section 4.2]{MWblob}).

\par
Given $q\in \mathbb{C}\backslash\{0\}$ and $n\in \mathbb{Z}_{> 0}$, $[n]_q\in \mathbb{C}$ denotes the quantum number. That is,
\begin{equation}
 [n]_q= q^n+q^{n-2} + \ldots + q^{-n+2} + q^{-n}.
 \end{equation}

\begin{defi} \label{defi blob algebra}
Let $q \in \mathbb{C}\backslash\{0\}$ and $m \in \mathbb{Z}_{> 0}$ such that $[m]_q \neq 0$. The blob algebra $b_n(q,m)$ is the $\mathbb{C}$-algebra associated with the generators $1$, $U_0$, $U_1$, \ldots , $U_{n-1}$ and the relations
\begin{align}
  U_i^2 &= -[2]_qU_i,  &   &   \mbox{if } 1 \leq i \leq n-1; \\
  U_iU_jU_i&= U_i,   &   &   \mbox{if } |i-j|=1 \mbox{ and } i,j\neq 0;\\
  U_iU_j& = U_jU_i,  &   &   \mbox{if } |i-j|>1;\\
  U_1U_0U_1&=-[m-1]_qU_1; & &   \\
  U_0^{2} &= -[m]_q U_0; & &
\end{align}
\end{defi}

To establish a connection between the classical presentation for the blob and the graded algebra, some restrictions should be imposed on the parameters $q$ and $m$.

\begin{condi}     \label{conditions}
Henceforth, an $l$-th primitive root of unity $q\in \mathbb{C}$ and an integer $m $ satisfying  $1\leq m <l$ are fixed. Furthermore, it is assumed that $l$ is odd and
\begin{equation}  \label{equ conditions q and m}
  q^m\neq q^{-m},   \qquad  q^m\neq q^{-m+2},   \qquad q^{-m}\neq q^{m+2}.
\end{equation}
\end{condi}

It should be noted that each of the above conditions can be considered a congruence. Specifically, (\ref{equ conditions q and m}) can be restated as
\begin{equation}  \label{condition second equation}
 m \not \equiv 0 \mod l,    \qquad  m \not \equiv  1 \mod l,   \qquad m \not \equiv - 1 \mod l .
\end{equation}

 Let $I=\mathbb{Z}/ l \mathbb{Z}$. To describe the $\mathbb{Z}$-grading on the graded presentation of $b_n(q,m)$, it is convenient to introduce the matrix
$(a_{ij})_{i,j \in I}$, given by
$$ a_{ij}= \left\{
      \begin{array}{rl}
        2, & \hbox{if } i=j \mod l; \\
        0, & \hbox{if } i \not = j \pm 1 \mod l ; \\
        -1, & \hbox{if } i = j\pm 1 \mod l.
      \end{array}
    \right.
  $$

 \begin{teo}  \cite[Corollary 3.6]{PR-Hblob}\label{thm isomorphism blob KRL}
The blob algebra $b_n(q,m)$ is isomorphic to the $\mathbb{C}$-algebra generated by
$$ \{ \psi_1,\cdots , \psi_{n-1} \} \cup  \{ y_1,\cdots , y_{n} \}    \cup \{ e(\textbf{i}) \mbox{         }  | \mbox{         } \textbf{i} \in I^n \} $$
subject to the following relations for $\textbf{i},\textbf{j}\in I^n$ and all admissible $r, s$
\begin{align}
\label{aaaa}       y_1e(\textbf{i})  & = 0  & \mbox{ if } i_1 = \pm k \mod{l}  \\
\label{bbbb}       e(\textbf{i})    & = 0 &  \mbox{ if } i_1 \neq  \pm k \mod{l}  \\
\label{kl2} e(\textbf{i})e(\textbf{j})& =\delta_{\textbf{i,j}} e(\textbf{i}), & \\
\label{kl3}\sum_{\textbf{i} \in I^n} e(\textbf{i})& =1, \\
\label{kl4}y_{r}e(\textbf{i})& =e(\textbf{i})y_r ,& \\
\label{kl5}\psi_r e(\textbf{i})& =e(s_r \textbf{i}) \psi_r, \\
\label{kl6}y_ry_s& = y_sy_r,&
&   \\
\label{kl7}\psi_ry_s& = y_s\psi_r,&
 \mbox{ if } s\neq r,r+1 \\
\label{kl8}\psi_r\psi_s& = \psi_s\psi_r,&
 \mbox{ if } |s-r|>1 \\
\label{kl9}\psi_ry_{r+1}e(\textbf{i})& =\left\{ \begin{array}{l}
(y_r\psi_r +1)e(\textbf{i}) \\
y_r\psi_r e(\textbf{i})   \\
\end{array}
\right.  &\begin{array}{l}
\mbox{if   }  i_r=i_{r+1}   \\
\mbox{if   }   i_r \neq i_{r+1}     \\
\end{array} \\
\label{kl10}y_{r+1}\psi_re(\textbf{i})& =\left\{ \begin{array}{l}
(\psi_ry_r +1)e(\textbf{i}) \\
\psi_ry_r e(\textbf{i})   \\
\end{array}
\right. &  \begin{array}{lcc}
\mbox{if   }  i_r=i_{r+1}  \\
\mbox{if   }   i_r \neq i_{r+1}   \\
\end{array}
\end{align}
\begin{align}
\label{kl11}\psi_r^{2}e(\textbf{i})& =\left\{ \begin{array}{l}
0   \\
e(\textbf{i})  \\
(y_{r+1}-y_{r})e(\textbf{i}) \\
(y_{r}-y_{r+1})e(\textbf{i})  \\
\end{array}
\right. &  \begin{array}{l}
\mbox{if   }  i_r=i_{r+1}  \\
\mbox{if   }   i_r \neq i_{r+1} \pm 1    \\
\mbox{if   }  i_{r+1}= i_r+1   \\
\mbox{if   }  i_{r+1}= i_r-1  \\
\end{array} \\
\label{kl12}\psi_r\psi_{r+1}\psi_re(\textbf{i})& =\left\{ \begin{array}{l}
(\psi_{r+1}\psi_r\psi_{r +1} +1)e(\textbf{i})  \\
(\psi_{r+1}\psi_r\psi_{r +1} -1)e(\textbf{i})    \\
(\psi_{r+1}\psi_r\psi_{r +1} )e(\textbf{i})    \\
\end{array}
\right. &  \begin{array}{l}
\mbox{if   }  i_{r+2}=i_r=i_{r+1}-1   \\
\mbox{if   } i_{r+2}=i_r=i_{r+1}+1      \\
\mbox{otherwise   }        \\
\end{array}  \\
 \label{blob relation} e(\textbf{i})&=0  &  \mbox{if } \left\{
                                \begin{array}{l}
                                 i_1=k \mbox{ and } i_2=k-1, \mbox{ or }  \\
                                 i_1=-k \mbox{ and } i_2=-k-1.
                                \end{array}
                              \right.
\end{align}
where
\begin{itemize}
\item $k\in \mathbb{Z}$ satisfies  $2k \equiv  m \mod l$;
\item for a sequence $\textbf{i} \in I^n$, $i_j\in I $ denotes the $j$-th coordinate of $\textbf{i}$;
\item $ s_r := (r,r+1)$ is the simple transposition acting on $ I^n $ by permutation of the coordinates $ r, r+1$.
\end{itemize}
The conditions
$$   \begin{array}{ccccc}
       deg\mbox{  } e(\textbf{i})=0,  &  &  deg\mbox{  } y_r=2,  &  &
 deg\mbox{  }\psi_se(\textbf{i})=-a_{i_{s},i_{s+1}} \end{array}
$$
for $1 \leq r\leq n$,  $1 \leq s \leq n-1 $, and $\textbf{i} \in I^{n}$
define a unique $\mathbb{Z}$-grading on $ b_n(q,m) $.
 \end{teo}

\begin{rem}\rm
Brundan and Kleshchev constructed in \cite{brukle} isomorphisms between cyclotomic Hecke algebras and cyclotomic KLR algebras of type $A$. On the other hand, $b_n(q,m)$ is known to be a quotient of the cyclotomic Hecke algebra of type $G(2,1,n)$. These two facts were used in \cite{PR-Hblob} in order to obtain the graded presentation of $b_n(q,m)$. In this setting, the relevant KLR algebra can be recovered from $b_n(q,m)$ by drooping (\ref{blob relation}). 
\end{rem}

It turns out that working with a variation of the diagrammatic calculus introduced by Khovanov and Lauda \cite{KLdiagram} is easier than working with the presentation given in Theorem \ref{thm isomorphism blob KRL}. In this setting, long algebraic computations involving relations are reduced to a series of diagrams that most transparently encode the same information. However, it should be noted that the diagrammatic approach still has certain subtleties. For example, as will be seen later, given an arbitrary diagram, it is a difficult problem to decide whether it is zero or not.

\par
 A \emph{Khovanov--Lauda diagram}, or a diagram for short, consists of $n$ points on each of two parallel edges, the top edge and the bottom edge, that are connected by $n$ lines. Each line must connect a node on the top edge with a node on the bottom edge. Lines can intersect, but no triple intersection is allowed. Each line can also be decorated by dots, but dots cannot be located at the intersection of two lines. Finally, each diagram is labeled by a sequence $\bs{i}=(i_1, \ldots , i_n) \in I^n$.  This sequence is written under the bottom edge of the diagram. An example of such a diagram is shown below.

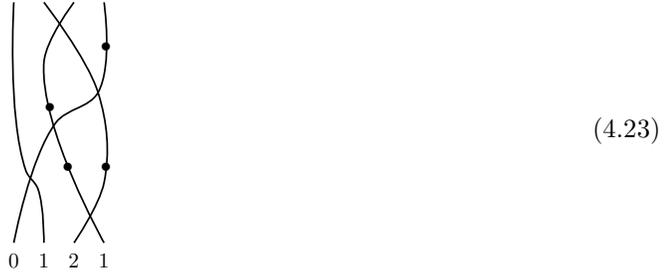
\begin{equation}  \label{ejemplo  KL diagrama}
\psscalebox{.8}{\wordtodiagram }
\end{equation}

Given a diagram $D$, $b(D)$ denotes its bottom sequence. Once $b(D)$ has been fixed, a top sequence for $D$, namely, $t(D)$, is automatically defined in an obvious manner. For instance, if $D$ is the diagram  (\ref{ejemplo  KL diagrama}), then
$t(D) =(1,2,1,0)$. $\mathbb{C}$-linear combinations of these diagrams are now considered modulo planar isotopy and the following relations.\footnote{Relations (\ref{first diagrammatic relation})--(\ref{last diagrammatic relation}) depend on Conditions \ref{conditions}. Moreover, it should be noted that for a fixed choice of the parameters $q$ and $m$, $k \in I$ is determined by the congruence $2k \equiv m \mod l$.}

\begin{equation}  \label{first diagrammatic relation}
 \yUnoIdempotenteCero =0
 \end{equation}

\begin{equation}
 \IdempotenteCeroHecke  =0 \qquad \mbox{ if }   \left\{
\begin{array}{l}
   i_1\neq \pm k ; \\
   i_1=k  \mbox{ and } i_2=k-1;      \\
  i_1=-k  \mbox{ and } i_2=-k-1.
\end{array}
     \right.
 \end{equation}

\begin{equation}  \label{puntocruzando de izquierda a derecha}
 \YconPsiUno  \quad = \YconPsiDos \quad + \delta_{i_r,i_{r+1}}\YconPsiTres
 \end{equation}

\begin{equation}   \label{puntocruzando de derecha a izquierda}
 \YconPsiUnoa \quad = \YconPsiDosa \quad + \delta_{i_r,i_{r+1}} \YconPsiTresa
 \end{equation}

\begin{equation}  \label{diagrammatic relation quadratic equal to zero}
\PsiCuadradoUno  =0 \mbox{ if } i_r=i_{r+1}   \end{equation}

 \begin{equation}   \label{diagrammatic relation quadratic equal to identity}
     \PsiCuadradoUno \quad  = \YconPsiTres   \quad \mbox{ if } i_r\neq i_{r+1} \pm 1 \mbox{ and }  i_r \neq i_{r+1}.
\end{equation}

  \begin{equation}   \label{diagrammatic relation quadratic two dots}
     \PsiCuadradoUno \quad  = \pm \left(  \PsiCuadradoDos \quad  -  \PsiCuadradotres \quad \right)  \quad \mbox{ if } i_{r+1}= i_r \mp 1
\end{equation}

  \begin{equation}  \label{last diagrammatic relation}
\TrenzaUno\quad  = \quad \TrenzaDos \quad   +   \alpha  \quad \TrenzaTres
\end{equation}
where $$\alpha = \left\{ \begin{array}{rl}
1, &\mbox{if   }  i_{r+2}=i_r=i_{r+1}-1;   \\
-1, & \mbox{if   } i_{r+2}=i_r=i_{r+1}+1;      \\
0, &\mbox{otherwise.  }        \\
\end{array}
\right.   $$

Finally, a multiplication operation on diagrams is defined by concatenation. More precisely, given two diagrams $D_1$ and $D_2$, their product is defined by
\begin{equation}
D_1D_2 = \left\{
\begin{array}{rl}
 0, &  \mbox{ if } b(D_1) \neq t(D_2); \\
  D_1 \circ D_2,  & \mbox{ if }     b(D_1) = t(D_2),
\end{array}   \right.
\end{equation}
where $D_1 \circ D_2$ is the diagram obtained from $D_1$ and $D_2$ by identifying the bottom of $D_1$ with the top of $D_2$. This multiplication is extended to $\mathbb{C}$-linear combinations of diagrams by linearity. The connection between $b_n(q,m)$ and diagrams is given by the following assignment:
\begin{equation}  \label{assignment}
 \KLRidempotent  \qquad   \qquad \KLRY   \qquad   \qquad  \KLRpsi
\end{equation}

By using this assignment, a $\mathbb{C}$-linear combination of diagrams can be associated with each element in $b_n(q,m)$ and vice versa. Furthermore, by the relations imposed on the diagrams, equivalent diagrams are associated with the same element in $b_n(q,m)$.

\section{ The intermediate sequence principle.}

The purpose of this section is to provide a sufficiency criterion for determining whether a diagram is zero. It will be called the \emph{Intermediate Sequence Principle}.  To this end, the construction of a graded cellular basis for $b_n(q,m)$ given in \cite{PR-Hblob} should be recalled. To properly introduce this basis, certain combinatorial objects should be defined, namely, one-line bipartitions and bitableaux, as well as residue sequences. An alternative interpretation of bitableaux as walks on the Pascal's triangle is also introduced. This section is concluded by showing with an example how the aforementioned criterion can be used to simplify arguments involving diagrams.

\subsection{Combinatorics}

Throughout this section, a positive integer $n$ is fixed. A \emph{one-line bipartition} of $n$ is an ordered pair $\lambda =(\lambda_1, \lambda_2)$, where $\lambda_1$ and $\lambda_2$ are non-negative integers such that $\lambda_1 + \lambda_2 =n$. The set of all one-line bipartitions of $n$ is denoted by $\bip $. It is useful to identify a one-line bipartition $\lambda$ with its Young diagram $[\lambda ]$, defined as
\begin{equation}
[\lambda ]= \{ (1,j,k) \in \{ 1\} \times \mathbb{N} \times \{1,2 \} \mbox{ } |  \mbox{ } 1\leq j \leq \lambda_k \}.
\end{equation}
$[\lambda]$ can be visualized as an ordered pair of usual one-line Young diagrams. A $\lambda$-\emph{bitableau} is a bijection
$\mathfrak{t}: [\lambda] \rightarrow \{1, \ldots , n\} $. $\mathfrak{t} $ is said to have shape $\lambda$, and the notation $\operatorname{Shape}(\mathfrak{t})=\lambda$ is used. $\mathfrak{t} $ may be regarded as a labeling of the boxes in $[\lambda ]$ using elements from $\{1,\ldots ,n\}$. A $\lambda$-bitableau is called standard if in each component its entries increase from left to right. The set of all standard $\lambda $-bitableaux is denoted by $\operatorname{Std}(\lambda)$. Moreover, let $\operatorname{Std}(n):= \bigcup_{\lambda \in \bip}\operatorname{Std}(\lambda)$. Examples of standard $\lambda $-bitableaux are shown in (\ref{ejemplo bitableau mayor}) and (\ref{ejemplo bitableau}), where $n=9$ and $\lambda =(6,3)$. Given $1\leq k \leq n$ and $\mathfrak{t} \in \operatorname{Std}(\lambda)$, $\mathfrak{t}_{|_k}$ denotes the standard bitableau obtained from $\mathfrak{t}$ by erasing the boxes in $\mathfrak{t}$ with entries greater than $k$.

\begin{equation}  \label{ejemplo bitableau mayor}
\left( \bitableauone \mbox{ , }\bitableaudos \right)
\end{equation}

\begin{equation}   \label{ejemplo bitableau}
\left( \bitableautres \mbox{ , }\bitableaucuatro \right)
\end{equation}

For $\lambda = (\lambda_1, \lambda_2) \in \bip $, a bitableau of particular interest is $\mathfrak{t}^{\lambda} \in \operatorname{Std}(\lambda)$, which is defined as follows. Let $\mu = \min \{\lambda_1 , \lambda_2\}$. Then, the numbers $1,2,\ldots , n$  are located in $[\lambda]$ and are in increasing order along the rows according to the following rules:\\\\
- Odd numbers less than $2\mu$ are located in the second component of $\mathfrak{t^{\lambda}}$.\\
- Even numbers less than or equal to $2\mu$ are located in the first component of $\mathfrak{t^{\lambda}}$.\\
- Numbers greater than $2\mu$ are located in the remaining boxes.

\medskip
For example, if  $n=9$ and $\lambda =(6,3)$, then $\mathfrak{t}^{\lambda}$ corresponds to the bitableaux shown in (\ref{ejemplo bitableau mayor}).

\medskip
The symmetric group  $\mathfrak{S}_n$ acts naturally (on the left) on the set of $\lambda$-bitableaux by permuting the entries. Given $\lambda \in \bip $ and $\mathfrak{t} \in \operatorname{Std}(\lambda)$,  $d(\mathfrak{t})\in \mathfrak{S}_n$ denotes the permutation that satisfies $ \mathfrak{t}=  d(\mathfrak{t})\mathfrak{t}^{\lambda}  $. For instance, if $\mathfrak{t}$ denotes the bitableau in (\ref{ejemplo bitableau}), then
$  d(\mathfrak{t}) = s_6s_5s_7s_4s_6s_8s_3s_5s_7s_2s_4s_6s_1s_3s_5 $.

\begin{rem}  \rm  \label{remark independence of reduced expression}
By the definition of $\mathfrak{t}^{\lambda}$, the permutation $d(\mathfrak{t})$ is $321$-avoiding for all $\mathfrak{t} \in \operatorname{Std}(\lambda)$. That is, any  two reduced expressions of $d(\mathfrak{t})$ are related through a series of Coxeter relations of type $s_is_j =s_js_i $ for $|i-j|>1$.
\end{rem}

\par
A sequence $\bs{i}^{\mathfrak{t}} \in I^n$ will be associated with each $\mathfrak{t} \in \operatorname{Std}(n)$ and will be called the \emph{residue sequence} of $\mathfrak{t}$. It is first recalled that a primitive $l$-th root of unity $q\in \mathbb{C}$ and a positive integer $m$ satisfying Conditions \ref{conditions} have been fixed.

\begin{defi} \label{defi sequence residues}
Let $\lambda \in \bip$, and let $A=(1,c,d)$ be a box in $[\lambda]$. It is assumed that $k \in \mathbb{Z} $ satisfies $2k \equiv m \mod l$. Then, the residue of $A$, $\rho (A) \in I$, is given by
\begin{equation}
\rho (A):= \left\{  \begin{array}{rl}
k+ (c-1),     &  \mbox{ if } d=1  ;    \\
-k + (c-1) ,  &    \mbox{ if } d=2.
\end{array}
  \right.
\end{equation}
Let $j$ be an integer with $1\leq j \leq n$.  For $\mathfrak{t} \in \operatorname{Std}(\lambda)$, the residue of $\mathfrak{t} $  at $j$ is defined by
$\rho_{\mathfrak{t}}(j)= \rho (A)$, where $A$ is the node occupied by $j$ in $\mathfrak{t} $.
Finally,  the residue sequence of $\mathfrak{t}$ is defined by $\bs{i}^{\mathfrak{t}}:= (\rho_{\mathfrak{t}}(1), \ldots , \rho_{\mathfrak{t}}(n)) \in I^n$. In particular, $\bs{i}^{\mathfrak{t}^\lambda }$ is simply written as $\bs{i}^{\lambda}$.
\end{defi}

\begin{exa}
If  $l=5$, $m=2$, and $\lambda = (6,3)$, then
\begin{equation}
\bs{i}^{\lambda}=\bs{i}^{\mathfrak{t}^{\lambda}}=(4,1,0,2,1,3,4,0,1)\in I^9.
\end{equation}
\end{exa}

It will now be shown that standard bitableaux may be regarded as walks on the Pascal triangle. This alternative interpretation allows providing two useful criterions for determining when two standard bitableuax give rise to the same residue sequence and when a standard bitableaux remains standard if a simple transposition acts on it.

\par
 The vertices on the Pascal triangle are labeled by ordered pairs of integers corresponding to level and column. For instance, (\ref{pascal triangle uno}) shows the Pascal triangle truncated at level $8$. The top vertex is labeled by $(0,0)$ and $\bullet = (5,-3)$.

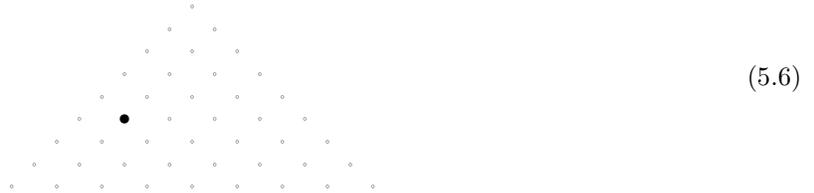
\begin{equation}  \label{pascal triangle uno}
\psscalebox{.3}{\pascalnueve }
\end{equation}

 Then, an $n$-walk on the Pascal triangle (or an $n$-walk for short) is defined to be a sequence of vertices $w =(v_0,v_1, \ldots , v_n)$ such that $v_{i} =(i, w_i)$ for all $0 \leq i \leq n$, $w_0=0$, and $w_{i}=w_{i-1} \pm 1$. Graphically, an $n$-walk may be regarded as a directed graph with vertices  $\{v_i\}_{0\leq i \leq n}$    and edges   $\{(v_i, v_{{i+1}})\}_{0\leq i \leq n-1}$. An $n$-walk $w^{\mathfrak{t}}=(v_0,v_1, \ldots , v_n)$ is associated with  $\mathfrak{t} \in \operatorname{Std}(n)$ as follows: $v_0=(0,w_0)=(0,0)$ and $v_i= (i, w_{i-1}+1)$  (resp.,  $v_i= (i, w_{i-1}-1)$) if $i$ is located in the first (resp., second) component of $\mathfrak{t}$. The $9$-walks associated with the  standard bitableaux in (\ref{ejemplo bitableau mayor}) and (\ref{ejemplo bitableau}) are drawn in (\ref{pascal triangle dos}) by a normal line and  a dashed line, respectively. It is an easy exercise to show that the above defines a bijection between $\operatorname{Std}(\lambda) $ and the set of all $n$-walks with final vertex $(n, \lambda_1 - \lambda_2)$.

 \begin{equation}  \label{pascal triangle dos}
\psscalebox{.3}{\pascalnueveA }
\end{equation}

\par
In view of Conditions \ref{conditions}, on the Pascal triangle, vertical lines are drawn on each column $j$ such that  $j \equiv -m \mod l$  (see Example \ref{example caminos similares al mayor} below). These lines will be called \emph{walls}.  Let $w$ be a $n$-walk that, at some vertex $v_i=(i, w_i)$, passes through a wall drawn on $j$. Then, let $w (i, j)$ be the $n$-walk obtained from $w $ by applying the reflection in the  wall drawn on $j$ to all vertices in $w $ after the vertex $v_i$. The notation $w \sim w'$ implies that there exists a sequence of walks
$(w^0, w^1 , \ldots , w^r)$ such that $w^0=w $, $w^r = w'$, and $ w^{s+1}=  w^s(i_s,j_s)$ for all $0\leq s <r$ and some integers $i_s$ and $j_s$.

\begin{exa}  \label{example caminos similares al mayor}  \rm
The Pascal triangle truncated at level $19$ is shown in  (\ref{dibujo caminos similares al mayor}). It is assumed that $l=5$ and $m=2$. The corresponding walls are drawn using thick vertical lines. The walks drawn correspond to all $19$-walks $w$ such that $w\sim w^\lambda$, where $w^\lambda$ is the $19$-walk associated with $\mathfrak{t}^{\lambda}$ for $\lambda =(2,17)$.
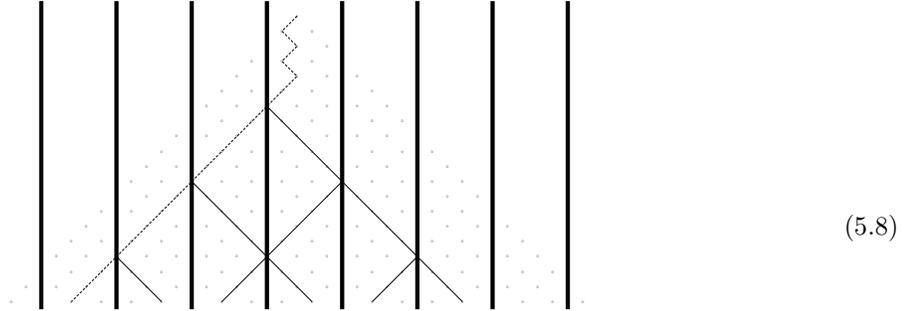
\begin{equation}  \label{dibujo caminos similares al mayor}
\psscalebox{.2}{\pascalB  }
\end{equation}
\end{exa}

\begin{lem} \label{lemma residues caminos}
Let $\mathfrak{s}, \mathfrak{t} \in \operatorname{Std}(n)$. Then, $\bs{i}^{\mathfrak{s}}=\bs{i}^{\mathfrak{t}}$ if and only if $w^{\mathfrak{s}} \sim w^{\mathfrak{t}}$.
\end{lem}

\begin{dem}
We begin by noticing that for any $\mathfrak{s} \in \operatorname{Std}(n)$ and $1\leq i \leq n $, we have
\begin{equation}  \label{resudios caminos}
2\rho_{\mathfrak{s}}(i) \equiv i-2 + (w^{\mathfrak{s}}_i -w^{\mathfrak{s}}_{i-1}) (w^{\mathfrak{s}}_i +m) \mod l ,
\end{equation}
where $w^{\mathfrak{s}}=(v_0, \ldots, v_n)$ with $v_i= (i, w^{\mathfrak{s}}_i)$ is the $n$-walk associated with $\mathfrak{s}$. As $l$ is assumed to be odd, (\ref{resudios caminos})  determines uniquely the value of $\rho_{\mathfrak{s}}(i)$.

Let $\mathfrak{s}, \mathfrak{t} \in \operatorname{Std}(n)$ with $w^{\mathfrak{s}} \sim w^{\mathfrak{t}}$. We are going to show that $\bs{i}^{\mathfrak{s}}=\bs{i}^{\mathfrak{t}}$, this is, $\rho_{\mathfrak{s}}(i) = 	\rho_{\mathfrak{t}}(i)$, for all $1\leq i \leq n$.  We can assume without loss of generality that $w^{\mathfrak{t}}= w^{\mathfrak{s}} (a,j) $, for some integers $a$ and $j$. In this case, since the walks $w^{\mathfrak{s}} $ and $w^{\mathfrak{t}}$ coincide till the $a$-th step, we obtain 
\begin{equation} \label{resudios caminos A}
	\rho_{\mathfrak{s}}(i) = 	\rho_{\mathfrak{t}}(i),
\end{equation} 
for all $1\leq i \leq a$. Furthermore, we have 
\begin{equation} \label{resudios caminos B}
	w_a:=w_a^{\mathfrak{s}} =w_a^{\mathfrak{t}}\equiv  -m \mod l
\end{equation}
 and
\begin{equation}  \label{resudios caminos C}
	w_i^{\mathfrak{s}} - w_{i-1}^{\mathfrak{s}} = - (w_i^{\mathfrak{t}} - w_{i-1}^{\mathfrak{t}}), 
\end{equation}
for all $a < i \leq n$.  By combining (\ref{resudios caminos B}) and (\ref{resudios caminos C}) we obtain that 
\begin{equation} \label{resudios caminos D}
	w_{a+b}^\mathfrak{s} -w_{a} = -(w_{a+b}^\mathfrak{t} -w_{a}),
\end{equation}
for all $b\geq 1$ such that $a+b \leq n$. Therefore, 
\begin{equation} \label{resudios caminos E}
	w_{a+b}^\mathfrak{s} + m \equiv -(w_{a+b}^\mathfrak{t} +m) \mod l.
\end{equation}
Then, a combination of (\ref{resudios caminos}), (\ref{resudios caminos C}) and (\ref{resudios caminos E}) yields $ \rho_{\mathfrak{s}}(a+b) = 	\rho_{\mathfrak{t}}(a+b)$. Thus, $ \bs{i}^{\mathfrak{s}}=\bs{i}^{\mathfrak{t}} $. 

Conversely, suppose that $ \bs{i}^{\mathfrak{s}}=\bs{i}^{\mathfrak{t}}$. Let $a+1$ be the first step where the walks $w^{\mathfrak{s}} $ and $w^{\mathfrak{t}}$ differ (of course, if such a step does not exist then $w^\mathfrak{s} = w^\mathfrak{t}$ and there is nothing to prove).  The fact that $ \bs{i}^{\mathfrak{s}}=\bs{i}^{\mathfrak{t}}$ together with (\ref{resudios caminos}) force to the vertex $ (a, w^{\mathfrak{s}}_a)=(a, w^{\mathfrak{t}}_a)$ to be on a wall. We now consider the standard bitableau $\mathfrak{u}$ defined by the condition:  $w^{\mathfrak{u}}:= w^{\mathfrak{s}}(a, w^{\mathfrak{s}}_a) $. If $w^\mathfrak{u}= w^\mathfrak{t}$ then we are done. Otherwise, we have $ \bs{i}^{\mathfrak{u}}= \bs{i}^{\mathfrak{s}}=\bs{i}^{\mathfrak{t}}$ and the first step where the walks $w^{\mathfrak{u}} $ and $w^{\mathfrak{t}}$ differ is strictly greater than $a+1$. We can continue in this way till we eventually reach $\mathfrak{t}$ by applying a sequence of reflections. Therefore,  $w^{\mathfrak{s}} \sim w^{\mathfrak{t}}$.
\end{dem}

\medskip
This section is concluded by providing a criterion (in terms of walks) for determining when a standard bitableau remains standard if a simple transposition acts on it.

\begin{defi}  \rm
Let $\lambda \in \bip $ and $\mathfrak{s}, \mathfrak{t} \in \operatorname{Std}(\lambda) $. Let $$w^{\mathfrak{s}}=((0,w^{\mathfrak{s}}_0), \ldots , (n, w^{\mathfrak{s}}_n) \mbox{ and } w^{\mathfrak{t}}=((0,w^{\mathfrak{t}}_0), \ldots , (n, w^{\mathfrak{t}}_n))$$ be the walks associated with $\mathfrak{s}$ and $\mathfrak{t}$, respectively. $w^{\mathfrak{s}}$ is said to\emph{ have a hook at position $k$ } if $ w^{\mathfrak{s}}_{k-1}= w^{\mathfrak{s}}_{k+1} =w^{\mathfrak{s}}_k \pm 1$. Furthermore, ${\mathfrak{t}}$ is said to be obtained from ${\mathfrak{s}}$ by making a hook at position $k$ if
$w^{\mathfrak{s}}_j=w^{\mathfrak{t}}_j $ for all $j\neq k$ and $|w^{\mathfrak{s}}_k-w^{\mathfrak{t}}_k|=2$.
 \end{defi}

\begin{lem}  \label{lemma hook}
Let $\lambda \in \bip $ and $\mathfrak{s} \in \operatorname{Std}(\lambda) $. Then,  $s_k \mathfrak{s}$ is standard if and only if $\mathfrak{s}$ has a hook at position $k$.
\end{lem}
\begin{dem}
The bitableau $s_k \mathfrak{s}$ is standard if and only if  $k$ and $k+1$ are located in different components of $\mathfrak{s}$. However, $k$ and $k+1$ are located in different components of $\mathfrak{s}$ if and only if $\mathfrak{s}$ has a hook at position $k$.
\end{dem}

\subsection{KLR-Idempotents.}

Elements $e(\bs{i})\in b_n(q,m)$ are referred to as KLR-idempotents. In general, by working only with the presentation of $b_n(q,m)$ given in Theorem \ref{thm isomorphism blob KRL}, it cannot be determined whether an arbitrary KLR-idempotent is zero or not. Fortunately, by using the graded cellular basis of $b_n(q,m)$ constructed in  \cite{PR-Hblob}, an effective criterion can be provided for this problem. This criterion will be sufficient for proving the \emph{Intermediate Sequence Principle}.

\begin{defi}
Let $\lambda \in \bip $ and $\mathfrak{s},\mathfrak{t} \in \operatorname{Std}(\lambda) $. Moreover, let $d(\mathfrak{s})=s_{i_1} \dots s_{i_a}$ and $d(\mathfrak{t})= s_{j_1} \ldots s_{j_b}$ be reduced expressions for $d(\mathfrak{s})$ and $d(\mathfrak{t})$, respectively. Then, $\psi_{\mathfrak{st}}^\lambda$ is defined by
\begin{equation}
\psi_{\mathfrak{st}}^\lambda := \psi_{i_1} \ldots \psi_{i_a} e (\bs{i}^\lambda) \psi_{j_b}  \ldots \psi_{j_1}   \in b_n(q,m).
\end{equation}
\end{defi}

By (\ref{kl8}) and Remark \ref{remark independence of reduced expression}, the element $\psi_{\mathfrak{st}}^\lambda $ is well defined. That is, it  does not depend on the choice of the reduced expressions for $d(\mathfrak{s})$ and $d(\mathfrak{t})$. The following properties are an easy consequence of  (\ref{kl2}), (\ref{kl5}) and the definition of $\psi_{\mathfrak{st}}^\lambda $.

\begin{equation}   \label{ortogonality property}
e(\bs{i})  \psi_{\mathfrak{st}}^\lambda   = \left\{
\begin{array}{rl}
  \psi_{\mathfrak{st}}^\lambda,     &  \mbox{ if } \bs{i} = \bs{i}^{\mathfrak{s}};   \\
 0 ,    &      \mbox{ if } \bs{i} \neq \bs{i}^{\mathfrak{s}}.
\end{array}
  \right.  \qquad \psi_{\mathfrak{st}}^\lambda  e(\bs{i})    = \left\{
\begin{array}{rl}
  \psi_{\mathfrak{st}}^\lambda,     &  \mbox{ if } \bs{i} = \bs{i}^{\mathfrak{t}};   \\
 0 ,    &      \mbox{ if } \bs{i} \neq \bs{i}^{\mathfrak{t}}.
\end{array}
  \right.
\end{equation}

\begin{teo} \label{teo graded cellular basis blob}
\cite[Theorem 6.10]{PR-Hblob} The set $\{\psi_{\mathfrak{st}}^\lambda \mbox{ } | \mbox{ } \lambda \in \bip \mbox{, }   \mathfrak{s},\mathfrak{t} \in \operatorname{Std}(\lambda) \}$ is a graded cellular  basis of $b_n(q,m)$.
\end{teo}

\begin{lem}   \label{lemma sequence principle}
Let $\bs{i} \in I^n$. Then, $e(\bs{i}) \neq 0 \in b_n(q,m)$ if and only if there exist $\lambda \in \bip$ and $\mathfrak{t} \in \operatorname{Std}(\lambda)$ such that $\bs{i}=\bs{i}^\mathfrak{t}$.
\end{lem}

\begin{dem}
This is a direct consequence of (\ref{kl2}), (\ref{ortogonality property}), and Theorem \ref{teo graded cellular basis blob}.
\end{dem}

\medskip
Lemma \ref{lemma sequence principle} will now be translated into the diagrammatic setting. Of course, Lemma \ref{lemma sequence principle} implies that a diagram of type
\begin{equation}
\KLRidempotentsequenceprinciple
\end{equation}
is different from zero if and only if  its bottom sequence can be obtained as a residue sequence for some standard bitableaux. However, more can be said about an arbitrary diagram. Let $D$ be a diagram. A line is drawn in $D$ between (and parallel to) the top and bottom edges that does not intersect any intersection in $D$ (see (\ref{diagrama con linea}) below). This line determines a sequence in $I^n$. Obviously, this sequence depends on the height at which the line is drawn. Sequences obtained in this manner will be called \emph{intermediate sequences} of $D$.

\begin{equation}  \label{diagrama con linea}
\psscalebox{.8}{\wordtodiagramwithline } \longrightarrow   (1,1,0,2)
\end{equation}

\begin{teo} (Intermediate Sequence Principle (ISP))
Let $D$ be a diagram. It is assumed that an intermediate sequence $\bs{i}$ of $D$ cannot be obtained as the residue sequence of some standard bitableau. Then, $D=0$.
\end{teo}
\begin{dem}
Let $D' \in b_n(q,m)$ be the element associated with $D$ by  (\ref{assignment}). Then, $D'$ belongs to the ideal of $b_n(q,m)$ generated by  $e(\bs{i})$. It follows from Lemma \ref{lemma sequence principle} that the above ideal is zero. Therefore, $D'=0$.
\end{dem}

\medskip
This section is concluded by proving a proposition that shows that the ISP simplifies the diagrammatic calculus. The following lemma is required.

\begin{lem}   \label{lemma impossible subsequences}
Let $\bs{i}= (i_1,\ldots , i_n) \in I^n$. It is assumed  that $i_r=i_{r+1}=i_{r+2}+1$ or $i_{r+1}=i_{r+2}=i_r-1$ for some $1\leq r \leq n-2$. Then, $\bs{i} \neq \bs{i}^\mathfrak{t}$ for all $\mathfrak{t} \in \operatorname{Std}(n)$.
\end{lem}
\begin{dem}
The result will be proved only for $i_r=i_{r+1}=i_{r+2}+1$. The result for  $i_{r+1}=i_{r+2}=i_r-1$ is proved similarly. It is assumed that there exists $\mathfrak{t} \in \operatorname{Std}(n)$ such that $\bs{i}^\mathfrak{t} = \bs{i}$. Then, 
\begin{equation} \label{conditions on a particular lemma}
\rho_{t}(r)=\rho_{t}(r+1)=\rho_{t}(r+2)+1.
\end{equation}
It is first claimed that $r$ and $r+1$ are located in different components of $\mathfrak{t}$. To see this, if $r$ and $r+1$ were located in the same component of $\mathfrak{t}$, then we would have
\begin{equation}
0 = \rho_{t}(r+1)- \rho_{t}(r) \equiv 1 \mod l,
\end{equation}
which is absurd. It is now assumed that  $r$ and $r+2$ are located in the same component of $\mathfrak{t}$. Then, 
\begin{equation}
-1= \rho_{t}(r+2)- \rho_{t}(r) \equiv 1 \mod l,
\end{equation}
\noindent
which is impossible, as $l$ is assumed to be odd. A similar conclusion is obtained if it is assumed that $r+1$ and $r+2$ are located in the same component of
$\mathfrak{t}$. Therefore, $r+2$ cannot be located in the components of $\mathfrak{t}$. This implies that $\mathfrak{t}$ does not exist, proving the lemma.
\end{dem}

\begin{pro}
Relation (\ref{last diagrammatic relation}), or equivalently (\ref{kl12}), can be simplified as follows:

  \begin{equation}  \label{braid relation simplified one}
\TrenzaUno\quad  =  \alpha  \quad \TrenzaTres \mbox{ , }  \alpha = \left\{ \begin{array}{rl}
 -1 ,&  i_r=i_{r+2}=i_{r+1}+1; \\
 0 ,&  i_r=i_{r+2}=i_{r+1}-1.
\end{array}
  \right.
\end{equation}

  \begin{equation}  \label{braid relation simplified two}
 \TrenzaDos \quad=     \beta  \quad \TrenzaTres  \mbox{ , }  \beta = \left\{ \begin{array}{rl}
 -1 ,&  i_r=i_{r+2}=i_{r+1}-1; \\
 0 ,&  i_r=i_{r+2}=i_{r+1}+1.
\end{array}
  \right.
\end{equation}
\end{pro}

\begin{dem}
Only (\ref{braid relation simplified one}) will be proved. The diagram on the left-hand side of (\ref{braid relation simplified one}) is considered. If $i_r=i_{r+2}=i_{r+1}-1$, then the intermediate sequence determined by the dashed line
\begin{equation}\label{braid relation with intermediate sequence}
  \TrenzaAAA
\end{equation}
satisfies the hypothesis in Lemma \ref{lemma impossible subsequences}. By applying the ISP and Lemma \ref{lemma impossible subsequences}, the proof in this case is completed. It is now assumed that $ i_r=i_{r+2}=i_{r+1}+1$. By (\ref{last diagrammatic relation}), we have
  \begin{equation}  \label{trenza con ISP}
\TrenzaUno\quad  = \quad \TrenzaBBB \quad  - \quad \TrenzaTres \quad .
\end{equation}
The intermediate sequence determined by the dashed line in the right-hand side of (\ref{trenza con ISP}) satisfies the hypothesis in Lemma \ref{lemma impossible subsequences}. Then, by the ISP, the diagram located on the left of the right-hand side of (\ref{trenza con ISP}) is zero, and (\ref{trenza con ISP}) becomes (\ref{braid relation simplified one}).
\end{dem}

\section{Commutative subalgebras of $b_n(q,m)$.}

In the previous section, a residue sequence $\bs{i}^{\mathfrak{t}}\in I^{n}$ was associated with  $\mathfrak{t}\in \operatorname{Std}(n)$. In particular, given $\lambda \in \bip$, $\bs{i}^{\lambda}\in I^{n}$ was defined as the residue sequence corresponding to $\mathfrak{t}^\lambda$. Thus, a KLR-idempotent $ e(\bs{i}^{\lambda})\in b_n(q,m) $ can be associated with each $\lambda \in \bip $. Then, for each $\lambda \in \bip$, the idempotent truncation $b_n(\lambda):= e(\bs{i}^{\lambda}) b_n(q,m) e(\bs{i}^{\lambda})$ of $b_n(q,m)$ can be defined. In Section 3, a subalgebra of $\al$ isomorphic to the Dot--Line algebra was obtained. A subalgebra of $b_n(\lambda)$ isomorphic to the Dot--Line algebra will now be determined. The commutativity of the Dot--Line algebra motivates the definition of the following subalgebra of  $b_n(\lambda)$.
\begin{defi}   \label{defi Yla}
  Let $n\in \mathbb{N}$ and $\lambda \in \bip$. Then, $\Yla$ is defined as the subalgebra of $b_n(q,m)$ generated by
  $e(\bs{i}^\lambda), y_1e(\bs{i}^\lambda ), \ldots , y_ne(\bs{i}^\lambda )$.
\end{defi}
It is claimed that $Y_n(\lambda)$ is isomorphic to the Dot--Line algebra on $q_\lambda +1$ strings, where $q_\lambda$ is an integer that depends on $q$, $m$, and $\lambda$ (see Definition \ref{technical conditions on lambda} below). The proof of this result is postponed to the next section. In this section, attention is restricted to determining a generating set for $\Yla$ and to proving that this set satisfies the same relations as those satisfied by the generators of the Dot--Line algebra.
\par
To relate $\Yla$ to the Dot--Line algebra, the first objective is to show that several generators of $\Yla$ are  zero.  Specifically, if $\lambda = (\lambda_1 , \lambda_2)$ and $\mu =\min\{\lambda_1, \lambda_2 \}$, then $y_re(\bs{i}^\lambda )=0$ for all $1\leq r \leq 2 \mu$. To prove this, the following two lemmas are required.

\begin{lem}   \label{lema cruce de lineas uno}
Let $a$ be a positive integer such that $2a+2\leq n$. For $1\leq s \leq n$, $j_s \in I$ is defined by
\begin{equation} \displaystyle
j_s=\left\{  \begin{array}{rl}
-k +(s-1)/2,  & \mbox{ if } s \mbox{ is odd and } s\leq 2a; \\
 k +(s-2)/2,  & \mbox{ if } s \mbox{ is even and } s\leq 2a;  \\
k+a,  &   \mbox{ if } s=2a+1 ; \\
 -k+a +1 &    \mbox{ if } s=2a+2; \\
 \text{any value} & \mbox{ if } 2a+2 < s.
\end{array}
       \right.
\end{equation}
If $\bs{j}:=(j_1,\ldots , j_{n})\in I^n$, then $e(\bs{j})=0$.
\end{lem}

\begin{dem}
By Lemma (\ref{lemma sequence principle}), it is enough to show that there does not exist a one-line standard bitableau $\mathfrak{t}$ such that $\bs{i}^{\mathfrak{t}}= \bs{j}$. It is assumed toward a contradiction that such a bitableau $\mathfrak{t}$ exists. It is not difficult to see that Conditions \ref{conditions}  imply $\mathfrak{t}_{|_{2a}} = \mathfrak{t^\lambda}_{|_{2a}}$ for any $\lambda =(\lambda_1 , \lambda_2) \in \bip $ such that $\min \{\lambda_1 , \lambda_2 \} \geq a$. We recall that $\mathfrak{t}^{\lambda}_{|_{2a}}$ (resp. $\mathfrak{t}_{|_{2a}}$) is the standard  bitableau obtained from $\mathfrak{t}^{\lambda}$  (resp. $\mathfrak{t}$) by erasing the boxes in $\mathfrak{t}^{\lambda}$ (resp. $\mathfrak{t}$) with entries greater than $2a$. In particular, we have $\operatorname{Shape}(\mathfrak{t}_{|_{2a}}) =(a,a)$. Hence, the two possible locations for $2a+1$ in $\mathfrak{t}$ are $N_1=(1,a+1,1)$ and $N_2=(1,a+1,2)$. The residues for $N_1$ and $N_2$ are $k+a$ and $-k+a$, respectively. As $j_{2a+1}= k+a$ and
\begin{equation}
k+a \not \equiv -k+a \mod l,
\end{equation}
by Conditions \ref{conditions},  it follows that $2a+1$ must be located in $N_1$. The two possible locations for $2a+2$ in $\mathfrak{t}$ are  $N_2$ and $N_3=(1,a+2,1)$.  The residue of  $N_3$ is $k+a+1$. Therefore, at least one of the congruences
$$  j_{2a+2} \equiv -k+a \mod l,\quad   \mbox{ or } \quad j_{2a+2} \equiv k+a+1 \mod l   $$
must be satisfied. However, it follows from Conditions \ref{conditions} that this is impossible because $j_{2a+2} = -k+a+1$. Consequently, there does not exist a one-line standard bitableau $\mathfrak{t}$ such that $\bs{i}^{\mathfrak{t}}= \bs{j}$.
\end{dem}

\begin{lem}   \label{lema cruce de lineas dos}
Let $a$ be a positive integer such that $2a+2\leq n$. For $1\leq s \leq n$, $j_s \in I$ is defined by
\begin{equation} \displaystyle
j_s=\left\{  \begin{array}{rl}
-k +(s-1)/2,  & \mbox{ if } s \mbox{ is odd and } s\leq 2a+1; \\
 k +(s-2)/2,    & \mbox{ if } s \mbox{ is even and } s\leq 2a;  \\
 k+a +1 &    \mbox{ if } s=2a+2; \\
 \text{any value} & \mbox{ if } 2a+2 < s \leq n.
\end{array}
       \right.
\end{equation}
If $\bs{j}=(j_1,\ldots , j_{n})\in I^n$, then $e(\bs{j})=0$.
\end{lem}

\begin{dem}
This is proved by the same argument as the one given in the  proof of Lemma \ref{lema cruce de lineas uno}.
\end{dem}

\begin{teo}   \label{teo yr zero uno}
Let  $\lambda=(\lambda_1 , \lambda_2) \in \bip$ and $\mu =\min\{\lambda_1, \lambda_2 \}$.  Then, $y_re(\bs{i}^\lambda )=0$ for all $1\leq r \leq 2 \mu$.
\end{teo}
\begin{dem}
We proceed by induction on $r$. It should first be noted that $y_1e(\bs{i}^\lambda )=0$ by (\ref{first diagrammatic relation}). Let $\bs{i}^{\lambda}=(i_1,\ldots , i_n) $. By the definition of $\mathfrak{t}^\lambda $ and its residue sequence, it is straightforward to see that
\begin{equation}   \label{residues yr equal to zero}
i_r=\left\{  \begin{array}{rl}
-k +\frac{r-1}{2},  & \mbox{ if } r \mbox{ is odd and } r\leq 2\mu; \\
 &   \\
 k +\frac{r-2}{2},    & \mbox{ if } r \mbox{ is even and } r\leq 2\mu,
 \end{array}
       \right.
\end{equation}
\noindent
for all $1 \leq r \leq 2\mu$. Therefore, Conditions \ref{conditions} imply $i_{r}\neq i_{r+1}\pm 1$ for all $ 1\leq r \leq 2\mu$ with $r$ odd. That is, if $i_{r} = i_{r+1}\pm 1$  for some $ 1\leq r \leq 2\mu$ with $r$ odd, then by (\ref{residues yr equal to zero}), we would have $m  \equiv \pm 1 \mod l$, which violates (\ref{condition second equation}). Likewise, we have $i_r \neq i_{r+1}$. Therefore, by (\ref{diagrammatic relation quadratic equal to identity}), it follows that
 \begin{equation}    \label{y2igualcero}
     \PsiCuadradoUno \quad  = \YconPsiTres   \qquad \mbox{ for all } 1\leq r \leq 2\mu \mbox{ with }  r \mbox{ odd.}
\end{equation}
 In particular, we have
\begin{equation*}
y_2e(\bs{i}^\lambda )=   \ydosigualacero \quad \stackrel{ (\ref{y2igualcero})}{=}   \quad \ydosigualaceroa  \quad \stackrel{ (\ref{puntocruzando de derecha a izquierda})}{=}  \quad \ydosigualacerob \stackrel{ (\ref{first diagrammatic relation})}{=} 0
\end{equation*}
So far we have proven the theorem for $r=1$ and $r=2$. It is now inductively assumed that the theorem has been proven for $r$ and $r+1$, where $r$ is odd and $r+1< 2\mu$. To complete the proof, it suffices  to prove the theorem for $r+2$ and $r+3$. The following diagram is considered, as well as its intermediate sequence determined by the dashed line.

\begin{equation}   \label{diagram yr equal to zero}
\InductiveStepA
\end{equation}

By (\ref{residues yr equal to zero}), it is not difficult to see that if $2a=r-1$, then the intermediate sequence appearing in (\ref{diagram yr equal to zero}) has the form of the sequence in Lemma \ref{lema cruce de lineas uno}. Therefore, by the ISP, the diagram in (\ref{diagram yr equal to zero}) is equal to zero. Then, by using the fact that $y_{r}e(\bs{i}^\lambda )=0$, we obtain
\begin{align}  \label{yrmasdos igual a zero}
  0  = & \InductiveStepAA \stackrel{(\ref{diagrammatic relation quadratic two dots}) }{=} \InductiveStepB  \quad - \quad \InductiveStepBB    \nonumber \\
    \stackrel{\stackrel{(\ref{y2igualcero})}{(\ref{puntocruzando de izquierda a derecha})} }{=} &  \InductiveStepC \quad  -  \quad \InductiveStepCC   =   y_{r+2}e(\bs{i}^\lambda ) -y_{r}e(\bs{i}^\lambda ) = y_{r+2}e(\bs{i}^\lambda ).
\end{align}

To prove that $y_{r+3}e(\bs{i}^\lambda )=0$, it should be noted that applying Lemma \ref{lema cruce de lineas dos} and the ISP to the diagram below yields
\begin{equation}
\psscalebox{.8}{\InductiveStepD } =0
\end{equation}
Then,  by using the same method as in (\ref{yrmasdos igual a zero})  and the fact that $y_{r+1}e(\bs{i}^\lambda )=0$, it follows that  $y_{r+3}e(\bs{i}^\lambda )=0$.
\end{dem}

\medskip
\noindent
The previous theorem shows that several generators of $ \ydl $ are zero. The next goal is to prove that among the non-zero generators of $\ydl $, there are several repetitions. To properly describe which generators coincide, the following notation is introduced.

\begin{defi} \label{definition technical conditions}
Let $\lambda=(\lambda_1 , \lambda_2) \in \bip $ and  $\mu =\min \{\lambda_1 , \lambda_2\}$. Then, an integer is defined as follows:
\begin{equation}
a_\lambda := \left\{
              \begin{array}{ll}
                2\mu +m, & \hbox{ if } \mu = \lambda_1 ; \\
                2\mu +l-m, & \hbox{ if } \mu = \lambda_2.
              \end{array}
            \right.
\end{equation}
Furthermore, the following  technical conditions are imposed on $\lambda$
\begin{equation}  \label{technical conditions on lambda}
   n > a_\lambda \mbox{ and }
  n \not \equiv a_\lambda  \mod l.
\end{equation}
Under these conditions, there exist unique $q_\lambda , r_{\lambda} \in \mathbb{Z}_{\geq 0}$ such that
\begin{equation}
n-a_\lambda = q_\lambda l +r_\lambda \mbox{, with } 0 < r_\lambda < l .
\end{equation}
Finally, given $1\leq j \leq q_\lambda +1$, let
\begin{equation}\label{defin j underlined}
 \laba{j} := a_\lambda +(j-1)l.
\end{equation}
\end{defi}

The integer $a_\lambda$ corresponds to the level where for the first time the walk associated with $\mathfrak{t}^\lambda$ touches a wall of the fundamental alcove. The conditions on $\lambda $ that appear in (\ref{technical conditions on lambda}) can  also be rephrased in terms of alcoves and walls. Indeed, $ n > a_\lambda $ is equivalent to $\lambda $ being outside the fundamental alcove. Likewise, $ n \not \equiv a_\lambda \mod l$ is equivalent to $\lambda $ not being on a wall. Henceforth, these conditions on $\lambda$ will be tacitly assumed. Finally, it should be noted that the integers  $\laba{j}$ correspond to the set of all integers $r$ satisfying $2\mu < r \leq n$ such that $r\equiv a_\lambda \mod l$. This is equivalent to the condition that the integers $\laba{j}$ correspond to the levels where the walk associated with $\mathfrak{t}^{\lambda}$ touches a wall.

\begin{teo}  \label{teo yr zero dos}
Let $\lambda =(\lambda_1 , \lambda_2) \in \bip$ and  $\mu =\min \{\lambda_1 , \lambda_2\}$. Let $a_\lambda, q_{\lambda}$ and $r_\lambda$ as in Definition \ref{definition technical conditions}. Then,

\begin{equation}
y_re(\bs{i}^\lambda) = \left\{
                         \begin{array}{rl}
                           0, & \hbox{ if } 2\mu <r \leq \laba{1} ; \\
                           y_{\laba{j}+1}e(\bs{i}^\lambda ), & \hbox{ if } \laba{j} < r\leq \laba{j+1}; \\
                           y_{\laba{q_\lambda +1 }+1 }e(\bs{i}^\lambda)  , & \hbox{ if } \laba{q_\lambda +1} < r \leq n.
                         \end{array}
                       \right.
\end{equation}
\end{teo}

\begin{dem}
Let $\bs{i}^{\lambda} =(i^\lambda_1,\ldots , i^\lambda_n)$. The following is first proved.
\begin{cla}  \label{claim}
Let $2 \mu < r \leq n$. If $e(s_r \bs{i}^{\lambda})\neq 0$, then $r \equiv a_\lambda \mod l$.
\end{cla}
\begin{dem}
It is assumed that $e(s_r \bs{i}^{\lambda}) \neq 0$. By Lemma \ref{lemma sequence principle}, this is equivalent to the existence of a standard bitableau $\mathfrak{s} \in \mbox{Std}(n)$ such that $\bs{i}^{\mathfrak{s}}= s_r\bs{i}^{\lambda}$. It will be proved that if such a tableau $\mathfrak{s}$ exists, then $r \equiv a_\lambda \mod l$.
Let $ \bs{i}^\mathfrak{s} =(i^\mathfrak{s}_{1}, \ldots , i^\mathfrak{s}_n) $. Then, by the definition of $\mathfrak{s}$, we have
\begin{equation}   \label{eq demo Claim A}
i^\mathfrak{s}_s = \left\lbrace  \begin{array}{rl}
i^\lambda_{r+1},  &  \mbox{ if } s=r;\\
i^\lambda_{r},  &  \mbox{ if } s=r+1;\\
i^\lambda_{s},  &  \mbox{ otherwise.}
\end{array}      \right.
\end{equation}
It follows from (\ref{eq demo Claim A}) and the definition of $\mathfrak{t}^\lambda $ that
\begin{equation}  \label{eq demo Claim B}
i^\mathfrak{s}_{r+1}= \left\lbrace \begin{array}{r}
-k+\mu-r-1 \mbox{, if } \mu = \lambda_1 ; \\
k+\mu-r-1  \mbox{, if } \mu = \lambda_2 ;
\end{array}  \right.             \quad
i^\mathfrak{s}_r= \left\lbrace \begin{array}{r}
-k+\mu-r\mbox{, if } \mu = \lambda_1 ; \\
k+\mu-r\mbox{, if } \mu = \lambda_2 .
\end{array}  \right.
\end{equation}
It is claimed that $r$ and $r+1$ are located in different components of $\mathfrak{s}$. Indeed, if $r$ and $r+1$ were located in the same component of $\mathfrak{s}$, then we would have
\begin{equation}  \label{eq demo Claim C}
i^\mathfrak{s}_{r+1}-i^\mathfrak{s}_{r} \equiv 1 \mod l.
\end{equation}
Moreover, by (\ref{eq demo Claim B}), we have
\begin{equation}  \label{eq demo Claim D}
i^\mathfrak{s}_{r+1}-i^\mathfrak{s}_{r} \equiv -1 \mod l.
\end{equation}
By combining (\ref{eq demo Claim C}) and (\ref{eq demo Claim D}), a contradiction is derived because $l$ is assumed to be odd. Therefore, $r$ and $r+1$ are located in different components of $\mathfrak{s}$. This is equivalent by Lemma \ref{lemma hook} to $\mathfrak{s}$ having a hook at position $r$.

\par
Let now  $\mathfrak{t} = s_r \mathfrak{s}$. By the previous paragraph, $\mathfrak{t}$ is standard and $\bs{i}^\mathfrak{t} = \bs{i}^\lambda$. Furthermore, $\mathfrak{t}$ has a hook at position $r$. Let $w^{\lambda}$ (resp. $w^\mathfrak{t}$) denote the $n$-walk associated with $\mathfrak{t}^{\lambda}$ (resp. $\mathfrak{t}$). It is straightforward to see that $w^{\lambda}$ can be described as the $n$-walk that first zigzags on and off the central vertical axis of the Pascal's triangle by its negative part and then finishes as a straight line on the vertex $(n,\lambda_1 -\lambda_2)$ (see Example \ref{example caminos similares al mayor}). As $\bs{i}^\mathfrak{t} = \bs{i}^\lambda$, Lemma \ref{lemma residues caminos} implies  that $w^{\mathfrak{t}}$ has a hook at position $r$ only if $1\leq r\leq 2\mu$, or $2\mu < r \leq n$  and $r\equiv a_\lambda \mod l$.  This completes the proof of the claim.
\end{dem}

\medskip
Let us return to the proof of the theorem. Applying the ISP  and  Claim \ref{claim} yields
\begin{equation}
0 = \quad \psscalebox{.7}{\PsiCuadradoUnoA}  \stackrel{ (\ref{diagrammatic relation quadratic two dots})}{=}  \psscalebox{.7}{\PsiCuadradotresA } \quad  - \psscalebox{.7}{\PsiCuadradoDosA }   = y_{r+1}e(\bs{i}^{\lambda}) - y_{r}e(\bs{i}^{\lambda}),
\end{equation}
whenever $2 \mu < r \leq n$  and $r \not \equiv a_\lambda \mod l$. This proves the theorem for all $2\mu +1 < r \leq n $.

\par
 It remains to prove that $y_{2\mu +1} e(\bs{i}^{\lambda})=0$. To this end, it is convenient to distinguish two cases: $\mu = \lambda_1$ or $\mu = \lambda_2$. It should be noted that by Theorem \ref{teo yr zero uno},  $y_{2\mu -1} e(\bs{i}^{\lambda})=0$ and $y_{2\mu } e(\bs{i}^{\lambda})=0$.

\par
If  $\mu = \lambda_1$, then by using the same argument as in the proof of Theorem \ref{teo yr zero uno}, it follows that $y_{2\mu +1} e(\bs{i}^{\lambda})=0$ because $y_{2\mu -1} e(\bs{i}^{\lambda})=0$. By contrast, if $\mu = \lambda_2$, then the walk associated with $\mathfrak{t}^\lambda$ does not have a hook at position $2\mu$. Hence, applying the same argument as in the proof of Claim \ref{claim} yields
\begin{equation}\label{blablablablabla}
   y_{2\mu +1} e(\bs{i}^{\lambda})-y_{2\mu} e(\bs{i}^{\lambda})=0.
\end{equation}
In this case, it follows from (\ref{blablablablabla})that  $y_{2\mu +1} e(\bs{i}^{\lambda})=0$, as $y_{2\mu} e(\bs{i}^{\lambda})=0$.
\end{dem}

\medskip
Theorem \ref{teo yr zero uno} and Theorem \ref{teo yr zero dos} imply that  to generate $\Yla $, it is sufficient to consider the set of generators
 \begin{equation}\label{generators that does not work}
\{ e(\bs{i}^{\lambda}) \} \cup \{y_{r} e(\bs{i}^{\lambda}) \mbox{ } | \mbox{ }r=\laba{j} +1, 1\leq j \leq q_\lambda +1  \}.
 \end{equation}

However, it is more convenient for the present purposes to consider another set of generators rather than this.
\begin{defi}
Let $\lambda \in \bip $. Then, $ Y_j$ is defined by
\begin{equation}
Y_j =Y^{j,\lambda}_{q,m}:= \left(y_{\laba{j}+1} - y_{\laba{j}}  \right) e(\bs{i}^{\lambda}) \in \Yla,
\end{equation}
for all $1\leq j \leq q_\lambda +1$.
\end{defi}
Theorem \ref{teo yr zero dos} implies that
\begin{equation}  \label{old generators as a sum of new generators}
y_{\laba{j}+1}e(\bs{i}^{\lambda}) = \sum_{i=1}^{j}Y_i,
\end{equation}
for all $1\leq j \leq q_\lambda +1 $. Therefore, it follows from (\ref{generators that does not work}) and (\ref{old generators as a sum of new generators}) that the set
\begin{equation}\label{real generators}
  \{e(\bs{i}^{\lambda})\} \cup \{Y_{j} \mbox{ }| \mbox{ }  1\leq j \leq q_\lambda +1 \}
\end{equation}
generates $\Yla$. This section is concluded by showing that the elements $Y_j$ satisfy the same relations as the those satisfied by the generators of the Dot--Line algebra.

\begin{teo}   \label{teo homo sobre}
There exists a surjective algebra homomorphism $\Phi : DL_{q_{\lambda }+1} \rightarrow\Yla$ determined on the generators by $\Phi (X_j) = Y_j$.
\end{teo}

To prove Theorem \ref{teo homo sobre}, several lemmas are required.

\begin{lem}  \label{lemma yr con cruce}
Let $\lambda =(\lambda_1 , \lambda_2) \in \bip $ and  $\mu = \min \{\lambda_1 , \lambda_2\}$. Moreover, let $\bs{i}^\lambda =(i_1,\ldots, i_n) \in I^n$. Then, the elements $Y_j\in \Yla$ can be represented by the following  diagram
\begin{equation}   \label{diagrama Yr con dobles cruces}
\psscalebox{.5}{ \Yconcruce }
\end{equation}
for all $1\leq j \leq q_\lambda +1$.
\end{lem}
\begin{dem}
It is first noted that by the definition of $\mathfrak{t}^\lambda $ and its residue sequence, we have
\begin{equation}\label{resta resi}
  i_{r+1}-i_{r}\equiv 1 \mod l,
\end{equation}
for all $2\mu < r \leq n$. Then, by (\ref{diagrammatic relation quadratic two dots}) and (\ref{resta resi}), it follows that
\begin{equation}  \label{eq lema y con cruce}
Y_j = \left(y_{\laba{j}+1} - y_{\laba{j}}  \right) e(\bs{i}^{\lambda})= \quad \psscalebox{.6}{\PsiCuadradoUnoYY} ,
\end{equation}
for all $1\leq j \leq q_\lambda +1$. A part of a diagram of the form $ \quad \psscalebox{.2}{\doblecruce } $
is called a \emph{double crossing}. It will be  proved that (\ref{diagrama Yr con dobles cruces}) is equivalent to the right-hand side of (\ref{eq lema y con cruce}). To this end, each double crossing appearing in (\ref{diagrama Yr con dobles cruces}) is \emph{disarmed} from left to right by applying relations (\ref{diagrammatic relation quadratic equal to zero})--(\ref{diagrammatic relation quadratic two dots}).

\par
From left to right, we first have $a_\lambda -2\mu -1$ double crossings that can be reduced to straight lines by (\ref{diagrammatic relation quadratic equal to identity}). If the entire diagram has a unique double crossing, the proof is complete. If not, the leftmost double crossing is of the form (\ref{diagrammatic relation quadratic two dots}). Therefore, the diagram splits into two diagrams: one has a dot located on the $\laba{1}$-th line, and the other has a dot on the $(\laba{1}+1)$-th line. The following double crossing is of the form (\ref{diagrammatic relation quadratic equal to zero}); thus, the diagram with a dot on the $\laba{1}$-th line is zero. The leftmost double crossing in the remaining diagram has the form
\begin{equation}
\psscalebox{.6}{\doblecruceconpunto } .
\end{equation}
(\ref{puntocruzando de izquierda a derecha}) and (\ref{diagrammatic relation quadratic equal to zero}) yield
\begin{equation}  \label{xxx}
\psscalebox{.6}{\doblecruceconpunto } = - \quad\psscalebox{.6}{\doblecruceconpuntodesarmado } .
\end{equation}

In the resulting diagram, the leftmost crossing does not belong to a  double crossing but forms a braid as in (\ref{braid relation simplified two}). Therefore, it can be reduced to three straight lines. It is noted that the minus signs appearing in (\ref{braid relation simplified two}) and (\ref{xxx}) cancel each other. The resulting diagram has now double crossings of the type (\ref{diagrammatic relation quadratic equal to identity}), thus returning to the initial situation. Hence, the proof follows by repeating the same argument $j-1$ times.
\end{dem}

\begin{exa}
The proof of Lemma \ref{lemma yr con cruce} is here illustrated with an example. Let $n=21$ and $\lambda= (18,3)\in \operatorname{Bip_1}(21)$. It is assumed that $l=5$ and $m=2$. Under these conditions, we have $a_\lambda = 9$, $q_\lambda = 2$, and $r_\lambda =2$. Furthermore,
$$\bs{i}^{\lambda} =(4,1,0,2,1,3,4,0,1,2,3,4,0,1,2,3,4,0,1,2,3).$$
In this setting, we have $Y_1 = (y_{10}-y_9)e(\bs{i}^{\lambda} )$, $Y_2 = (y_{15}-y_{14})e(\bs{i}^{\lambda} )$, and $Y_3 = (y_{20}-y_{19})e(\bs{i}^{\lambda} )$. It will be shown that $Y_3$ can be obtained as a diagram of the form (\ref{diagrama Yr con dobles cruces}).
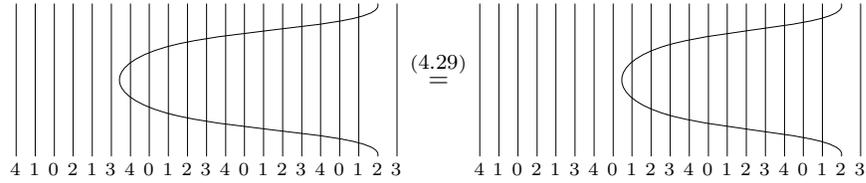
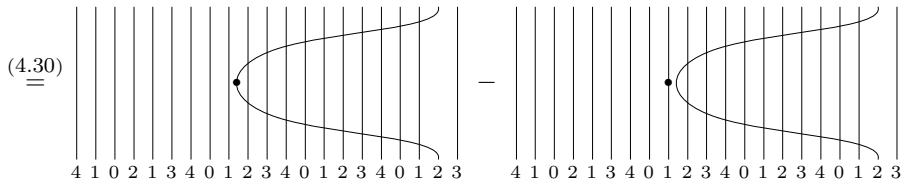
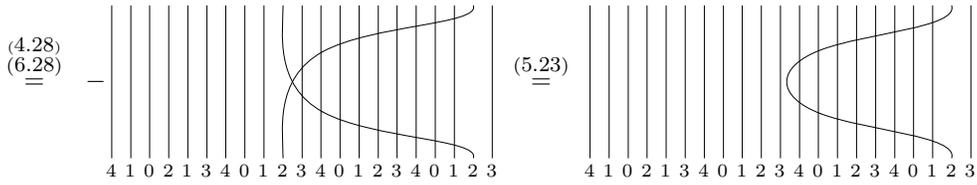
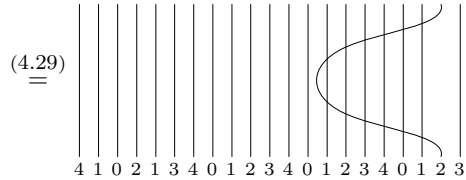
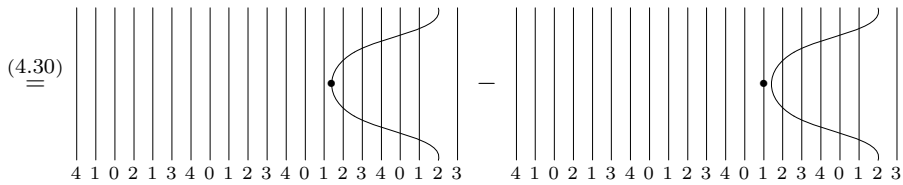
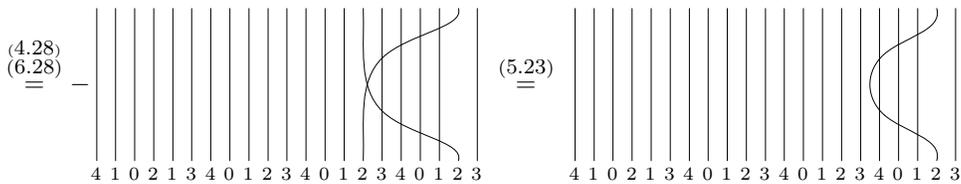
$$\hspace*{11mm}\dibuHA\stackrel{{\footnotesize(\ref{diagrammatic relation quadratic equal to identity})}}{=}\dibuHB$$
\\
$$\stackrel{{\footnotesize(\ref{diagrammatic relation quadratic two dots})}}{=}\dibuIA\;-\;\dibuIB$$
\\
$$\hspace*{2mm}\stackrel{\stackrel{{\footnotesize(\ref{diagrammatic relation quadratic equal to zero})}}{{\footnotesize(\ref{xxx})}}}{=}\hspace*{2mm}-\dibuJA\;\stackrel{{\footnotesize(\ref{braid relation simplified two})}}{=}\;\dibuJB$$
\\
$$\stackrel{{\footnotesize(\ref{diagrammatic relation quadratic equal to identity})}}{=}\dibuK$$
\\
$$\stackrel{{\footnotesize(\ref{diagrammatic relation quadratic two dots})}}{=}\dibuKA\;-\;\dibuKB$$
\\
$$\stackrel{\stackrel{{\footnotesize(\ref{diagrammatic relation quadratic equal to zero})}}{{\footnotesize(\ref{xxx})}}}{=}-\dibuLA\;\stackrel{{\footnotesize(\ref{braid relation simplified two})}}{=}\;\dibuLB$$
\\
$$\stackrel{{\footnotesize(\ref{diagrammatic relation quadratic equal to identity})}}{=}\dibuM=Y_3$$

\end{exa}

Lemma \ref{lemma yr con cruce} admits the following generalization.

\begin{lem}  \label{generalization lem con cruce}
Let $\lambda =(\lambda_1 , \lambda_2) \in \bip $ and  $\mu = \min \{\lambda_1 , \lambda_2\}$. For all $1\leq a\leq q_{\lambda}+1$, let $\hat{a}:=2\mu +(a-1)l$. Moreover, let $\bs{i}^\lambda =(i_1,\ldots, i_n) \in I^n$. Then, the elements $Y_j\in \Yla$ can be represented by the following  diagram
\begin{equation}   \label{diagrama Yr con dobles crucesgen}
\psscalebox{.5}{ \Yconcrucegen }
\end{equation}
for all $1\leq a \leq  j \leq q_\lambda +1$.
\end{lem}

\begin{dem}
The result follows by the same argument used in the proof of Lemma \ref{lemma yr con cruce}. The details are left to the reader.
\end{dem}

\begin{lem}  \label{lemma yr con cruce con punto zero}
Let $\lambda =(\lambda_1 , \lambda_2) \in \bip $, $\mu = \min \{\lambda_1 , \lambda_2\}$, and $\bs{i}^\lambda =(i_1,\ldots, i_n) \in I^n$. Then, the diagram
\begin{equation}   \label{diagrama Yr con dobles cruces A}
\psscalebox{.5}{ \YconcruceconpuntoZero }
\end{equation}
is equal to zero for all $1\leq j \leq q_\lambda +1$.
\end{lem}

\begin{dem}
Let $n'<n$. Given $\bs{a}=(a_1, \ldots , a_{n'}) \in I^{n'}$ and $\bs{b}=(b_1, \ldots , b_{n-n'})\in I ^{n-n'}$, let $\bs{a} \vee \bs{b} := (a_1, \ldots , a_{n'}, b_1 ,\ldots , b_{n-n'} )\in I^{n}$. By \cite[Section 4.4]{BKWGradedSpecht}, there exists an embedding $b_{n'}(q,m) \hookrightarrow b_n(q,m)$ determined  by
\begin{equation}
y_r \mapsto y_r, \quad \psi_s \mapsto \psi_s \mbox{ and } e(\bs{a})\mapsto \sum_{\bs{b} \in I^{n-n'}} e(\bs{a} \vee \bs{b})
\end{equation}
for all $1\leq r \leq n'$, $1\leq s < n'$, and $\bs{a} \in I^{n'}$. Thus, by the orthogonality of the KLR idempotents, if $y_re(\bs{a})=0$ in $b_{n'}(q,m)$, then $y_re(\bs{a}\vee\bs{b})=0$ in $ b_n(q,m) $ for all $\bs{b} \in I^{n-n'}$. This fact will be used to restrict attention to a subdiagram of (\ref{diagrama Yr con dobles cruces A}). The proof splits naturally into two cases:

\par
\emph{Case 1} $(\mu = \lambda_1)$. Let $\bs{a}:= (i_1, \ldots , i_{2\mu},i_{2\mu +1},i_{\laba{j}+1})\in I^{2\mu + 2}$. It will be proved that $\bs{a} = \bs{i}^{\nu }$, where $\nu = (\mu +1,\mu +1)$. To this end, it suffices to show that
\begin{equation} \label{eq residues for y con cruce con punto A}
i_{\laba{j}+1} - i_{2\mu } \equiv 1 \mod l.
\end{equation}
It is first noted that
\begin{equation}\label{eq residues for y con cruce con punto B}
i_{2\mu } \equiv k+\mu -1 \mod l.
\end{equation}
Moreover, as $\mu = \lambda_1$, we have
\begin{equation}\label{eq residues for y con cruce con punto C}
i_r\equiv -k+r-\mu -1 \mod l
\end{equation}
for all $2\mu < r \leq n$. In particular, 
\begin{equation}\label{eq residues for y con cruce con punto D}
  \begin{array}{rcll}
    i_{\laba{j}+1} & \equiv & -k+\laba{j}-\mu & \mod l, \\
                   & \equiv & -k+a_\lambda-\mu & \mod l, \\
                   & \equiv & -k+2\mu +m-\mu &\mod l, \\
                   & \equiv & k +\mu  & \mod l.
  \end{array}
\end{equation}
To obtain these congruences, $ 2k \equiv m \mod l$ and $\mu = \lambda_1$ were used. Then, (\ref{eq residues for y con cruce con punto A}) follows by subtracting (\ref{eq residues for y con cruce con punto B}) from (\ref{eq residues for y con cruce con punto D}). Thus,  $\bs{a} = \bs{i}^{\nu }$.

\par
 The subdiagram inside $\psscalebox{.4}{\rectPlomoClaro } $ in (\ref{eq y con cruce primer diagrama}) is now considered. By (\ref{eq residues for y con cruce con punto C}), we have $i_{2\mu +1} = -k+\mu $. It follows from (\ref{condition second equation}) and (\ref{eq residues for y con cruce con punto D}) that
 \begin{equation} \label{eq residues for y con cruce con punto E}
 i_{\laba{j}+1} - i_{2\mu +1} \equiv m  \not \equiv 0,1,-1 \mod l.
 \end{equation}
Hence, (\ref{puntocruzando de izquierda a derecha}) and (\ref{diagrammatic relation quadratic equal to identity}) can be applied to reduce the subdiagram inside $\psscalebox{.4}{\rectPlomoClaro } $ in (\ref{eq y con cruce primer diagrama}) to the element $y_{2\mu +2} e(\bs{a}) \in b_{2\mu +2} (q,m)$. Therefore, Theorem \ref{teo yr zero uno} implies $y_{2\mu +2} e(\bs{a})=0$, and by the first paragraph in this proof, it follows that the entire diagram (\ref{eq y con cruce primer diagrama}) is equal to zero.

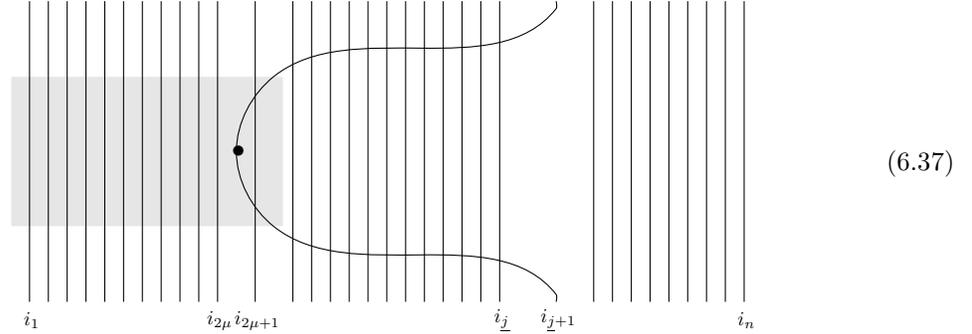
\begin{equation}   \label{eq y con cruce primer diagrama}
\psscalebox{.5 }{\YconcruceconpuntoZeroBox}
\end{equation}

\par
\emph{Case 2.} $(\mu= \lambda_2)$.  For this case, the subdiagram inside  $\psscalebox{.4}{\rectPlomoClaro } $ in (\ref{eq y con cruce segundo diagrama}) should be considered. This diagram corresponds to $y_{2\mu +1} e(\bs{a}) \in b_{2\mu +1} (q,m)$, where $\bs{a}:=(i_1,\ldots , i_{2\mu }, i_{\laba{j}+1} ) \in I^{2\mu +1}$. By repeating the same argument as in the previous case and by using Theorem \ref{teo yr zero dos}, it follows that the aforementioned subdiagram is equal to zero, proving the lemma in this case. The details are left to the reader.

\begin{equation}     \label{eq y con cruce segundo diagrama}
\psscalebox{.5 }{ \YconcruceconpuntoZeroBoxA }
\end{equation}

\end{dem}

\begin{lem}  \label{lemma y grande y chico}
Let $\lambda =(\lambda_1 , \lambda_2) \in \bip $ and  $\mu = \min \{\lambda_1 , \lambda_2\}$.  Then,

\begin{equation}  \label{Ygrandeychico}
Y_j y_{\underline{j} +1} = \left\lbrace
 \begin{array}{rll}
0,                         &   & \mbox{ if } j=1;    \\
\displaystyle -Y_j \sum_{i=1}^{j-1} Y_i, &   & \mbox{ if } 1< j \leq q_\lambda +1,
\end{array}
 \right.
\end{equation}
for all $1\leq j \leq q_\lambda +1$.
\end{lem}

\begin{dem}
Let $\bs{i}^\lambda =(i_1,\ldots, i_n) \in I^n$. By Lemma \ref{lemma yr con cruce}, the left-hand side of (\ref{Ygrandeychico}) can be represented diagrammatically by
\begin{equation}  \label{dibujoYgrandeychico}
Y_j y_{\underline{j} +1}= \psscalebox{.5}{\Ygrandeychico }
\end{equation}
The rectangle $\psscalebox{.4}{\rectPlomoClaro }$ in (\ref{dibujoYgrandeychico})  is drawn there to indicate a position in the diagram. The aim is to move the dot to $\psscalebox{.4}{\rectPlomoClaro  }$ by applying  (\ref{puntocruzando de izquierda a derecha}). It should be noted that the dot can freely move cross an intersection as long as the lines involved are labeled by different values; otherwise, the diagram splits into two diagrams. In this case, we have a diagram with a dot and another diagram with a braid crossing. For the diagram with the dot, (\ref{puntocruzando de izquierda a derecha}) can be applied again to move the dot to the left until it reaches a crossing formed by two lines labeled by the same value. At this point, the diagram splits into two diagrams again, one with a dot and another with a braid crossing. We can repeat this process to locate the dot in $\psscalebox{.4}{\rectPlomoClaro  }$. This diagram is equal to zero by Lemma \ref{lemma yr con cruce con punto zero}. Consequently, $Y_j y_{\underline{j} +1} $ is a linear combination of diagrams with a braid crossing (with the coefficients in the linear combination equal to one). By (\ref{braid relation simplified two}), each braid crossing can be disarmed into three straight lines multiplied by $-1$. Then, by Lemma \ref{lemma yr con cruce} and by using similar arguments to the those used in its proof, (\ref{Ygrandeychico}) is obtained.
\end{dem}

\begin{exa}
The previous lemma is now illustrated with an example. Let $n=21$ and $\lambda= (18,3)\in \operatorname{Bip_1}(21)$. It is assumed that $l=5$ and $m=2$.  It is first noted that
$$\hspace*{9mm}\dibuA \stackrel{{\footnotesize(\ref{puntocruzando de izquierda a derecha})}}{=}\dibuB$$
\\
$$\stackrel{{\footnotesize(\ref{puntocruzando de izquierda a derecha})}}{=}\dibuC \; + \; \dibuD$$
Both summands are now simplified as follows:
$$\dibuD \stackrel{{\footnotesize(\ref{braid relation simplified two})}}{=} - \dibuDA$$
\\
$$\stackrel{{\footnotesize\mbox{ Lemma } \ref{generalization lem con cruce}}}{=}-\dibuDB=-Y_3Y_2$$
Moreover,
$$\hspace*{9mm}\dibuC \stackrel{{\footnotesize(\ref{puntocruzando de izquierda a derecha})}}{=} \dibuCA$$
\\
$$\stackrel{{\footnotesize(\ref{puntocruzando de izquierda a derecha})}}{=}\dibuCC+\dibuCD$$
\\
$$\stackrel{\stackrel{{\footnotesize(\ref{puntocruzando de izquierda a derecha})}}{{\footnotesize(\ref{braid relation simplified two})}}}{=}\dibuCE-\dibuCF$$
\\
$$\stackrel{\stackrel{\footnotesize\mbox{ Lemma } \ref{lemma yr con cruce con punto zero}}{
{\footnotesize \mbox{ Lemma } \ref{generalization lem con cruce}}}}{=}-\dibuCG=-Y_3Y_1$$
Hence, 
$$\dibuA=-Y_3(Y_1+Y_2)$$
\end{exa}

\par
\begin{demespecial}
It should be shown that the elements $Y_j \in \Yla$ satisfy relations (\ref{dotlineone})--(\ref{dotlinethree}). $(\ref{dotlineone})$ follows directly from (\ref{kl4}), (\ref{kl6}), and the definition of the elements $Y_j$. Let now  $1<j \leq q_{\lambda}+1$. Then, we have

\begin{equation}\label{ysatisfiesdotlinerelations}
Y_j^2= Y_j(y_{\underline{j}+1}-y_{\underline{j}})e(\bs{i}^\lambda) = Y_j(y_{\underline{j}+1}-y_{\underline{j-1}+1})e(\bs{i}^\lambda)= -Y_j\sum_{i=1}^{j-1}Y_i - Y_j\sum_{i=1}^{j-1}Y_i =-2Y_j\sum_{i=1}^{j-1}Y_i,
\end{equation}
where the second equation is a consequence of Theorem \ref{teo yr zero dos} ( $y_{\underline{j}}e(\bs{i}^\lambda) = y_{\underline{j-1}+1}e(\bs{i}^\lambda)$), and  the third equation follows from Lemma \ref{lemma y grande y chico} and  (\ref{old generators as a sum of new generators}). Therefore, the elements $Y_j$ satisfy (\ref{dotlinetwo}).

\par
It remains to prove (\ref{dotlinethree}). To this end, it is first noted that  Theorem \ref{teo yr zero dos} yields $y_{\underline{1}} e(\bs{i}^\lambda)=0$. It follows now from Lemma \ref{lemma y grande y chico} that
\begin{equation}\label{ysatisfiesdotlinerelations dos}
  Y_1^2= Y_1y_{\underline{1}+1} e(\bs{i}^\lambda)-Y_1y_{\underline{1}} e(\bs{i}^\lambda)=0.
\end{equation} \end{demespecial}

\section{Isomorphism between $\Yla$ and Dot--Line algebra.}

In the previous section, it was shown that the generators $Y_j$ of $\Yla$ satisfy the same relations as those satisfied by the generators of the Dot--Line algebra. This allowed defining a surjective homomorphism $\Phi : DL_{q_{\lambda }+1} \rightarrow \Yla$ determined at the generators by $\Phi (X_j) = Y_j$. In this section, it will be shown that $\Phi $ is in fact an isomorphism.

\subsection{Dot--Line algebra as subalgebra of $b_n(q,m)$}
To establish the aforementioned isomorphism, it suffices to show that $\dim_\mathbb{C} \Yla = 2^{q_\lambda +1}$. The choice of a candidate for a basis of $\Yla$ should now be obvious. Namely,
\begin{equation}\label{basis of Y}
  \mathcal{Y}_n(\lambda):=\{ Y_1^{\alpha_1}\cdots  Y_{q_{\lambda }+1}^{\alpha_{q_{\lambda }+1}} \mbox{;  }  \alpha_i= 0,1 \text{, for all } 1\leq i \leq q_{\lambda }+1 \}.
\end{equation}
To prove that $\mathcal{Y}_n(\lambda)$ is a basis of $\Yla $, certain notations should be introduced.
\begin{defi}   \label{defi central reflection n walk}
Given $\lambda\in \bip$, the central reflection $n$-walk of $\lambda$ is denoted by $\omega_{\lambda}$ and is defined as the unique $n$-walk that satisfies the following two conditions:
\begin{enumerate}
  \item $\omega_{\lambda}\sim \omega^{\lambda}$.
  \item All nodes in $\omega_{\lambda}$ belong to the fundamental alcove or its walls.
\end{enumerate}
\end{defi}
$\mathfrak{t}_{\lambda}$ denotes the standard bitableaux associated with $\omega_{\lambda}$. In this setting, the first condition in Definition \ref{defi central reflection n walk} can be restated as $\bs{i}^{\mathfrak{t}_\lambda} = \bs{i}^\lambda$. It should be noted that $\operatorname{Shape}(\mathfrak{t}_{\lambda})\neq \lambda$, for all $\lambda \in \bip $. Henceforth, $\mu_\lambda \in \bip$ denotes the shape of the standard bitableau $\mathfrak{t}_{\lambda}$.

\begin{exa}\label{triej}
Definition \ref{defi central reflection n walk} will now be illustrated with an example.  Let $\lambda=(23,2)$, $l=7$, and $m=3$. The central reflection $25$-walk of $\lambda $ was drawn as a continuous black line. The $25$-walks drawn as gray lines and dashed lines correspond to the walks associated with  $\mathfrak{t}^{\lambda}$ and $\mathfrak{t}^{\mu_{\lambda}}$, respectively.
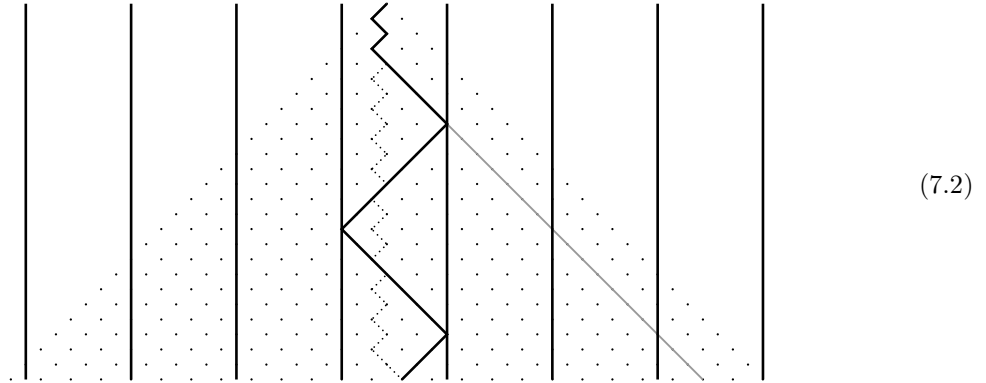
\begin{equation}\triejemB\end{equation}
\end{exa}

The interest in the standard bitableau $\mathfrak{t}_\lambda$ is justified by the following lemma.

\begin{lem}  \label{lema central reflection equal to product of Y's}
Let $\lambda \in \bip $. Then,
\begin{equation}\label{equ central reflection equal to product of Y's}
\psi^{\mu_{\lambda}}_{\mathfrak{t}_\lambda \mathfrak{t}_\lambda} =  \prod_{i=1}^{q_{\lambda}+1} Y_i.
\end{equation}
In particular, the right-hand side of (\ref{equ central reflection equal to product of Y's}) is an element of the graded cellular basis of $b_n(q,m)$, and therefore, different from zero.
\end{lem}

The proof of Lemma \ref{lema central reflection equal to product of Y's} is postponed to the next subsection.

\begin{teo}
  Let $\lambda \in \bip $. The set $ \mathcal{Y}_n(\lambda)$ is a basis of $\Yla$. In particular, $\dim_\mathbb{C} \Yla= 2^{q_{\lambda} +1} $. Consequently,  $\Phi : DL_{q_{\lambda }+1} \rightarrow \Yla $  is an isomorphism.
\end{teo}

\begin{dem}
As the generators $Y_j$ satisfy the same relations as those satisfied by the generators of the Dot--Line algebra, it follows that  $ \mathcal{Y}_n(\lambda)$ spans $Y_{n}^{\lambda}(q,m)$ by the same argument as that used in the proof of Lemma \ref{lema span set for Dot--Line}. Likewise, it follows that $ \mathcal{Y}_n(\lambda)$  is a linearly independent set by applying the same argument as that used in the proof of Theorem \ref{teo base de Dot line}. In this setting, the diagram (\ref{dibujo double dots in DL is diferent from zero}) is replaced by $\prod_{i=1}^{q_{\lambda}+1} Y_i$, which is different from zero by Lemma \ref{lema central reflection equal to product of Y's}.
\end{dem}

\subsection{Proof of Lemma \ref{lema central reflection equal to product of Y's}}

Let $\lambda \in \bip$. To prove Lemma \ref{lema central reflection equal to product of Y's}, it is required to have control over  $\psi^{\mu_{\lambda}}_{\T_\lambda \T_\lambda }$. By definition of this element, it suffices to have a simple process for obtaining a reduced expression for $d(\T_\lambda) \in \mathfrak{S}_n$. More generally, an algorithm for obtaining $d(\T ) \in \mathfrak{S}_n $ for all $\mathfrak{t} \in \operatorname{Std}(n)$ is reviewed. This process is described in
\cite[Section 4]{PR-Hblob}.

\begin{algo}  \label{algorithm}\mbox{}
Let $\mu \in \bip $ and $\T \in \operatorname{Std}(\mu)$. Let  $\omega^\mathfrak{t}$ and $\omega^\mu $ denote the $n$-walks associated  to $\T$ and $\T^\mu $, respectively. Then, the algorithm is as follows.
\begin{enumerate}
\item Draw both $\omega^\mathfrak{t}$ and $\omega_0:=\omega^\mu $ on the Pascal triangle. Set $d_0(\T)=1 $.
\item Assume that $\omega_{i-1}$ has been defined. Then, define $\omega_{i}$ as any $n$-walk obtained from $\omega_{i-1}$ by making a hook at some position $k$, where the area of the region bounded by $\omega^\mathfrak{t}$ and $\omega_{i}$ is lower than the area of the region bounded by $\omega_\mathfrak{t}$  and $\omega_{i-1}$. Set  $d_i(\T)=s_kd_{i-1}(\T
    )$.
\item Repeat $j$ times step $2$ until $\omega_j=\omega^\mathfrak{t}$. Then, $d(\T )=d_j(\T )$.
\end{enumerate}
\end{algo}

\begin{exa}
The standard bitableau $\T=\left( \;\young(57)\,,\,\young(12346)\;\right) $ is considered. Algorithm \ref{algorithm} is depicted in (\ref{eq algorithm}).
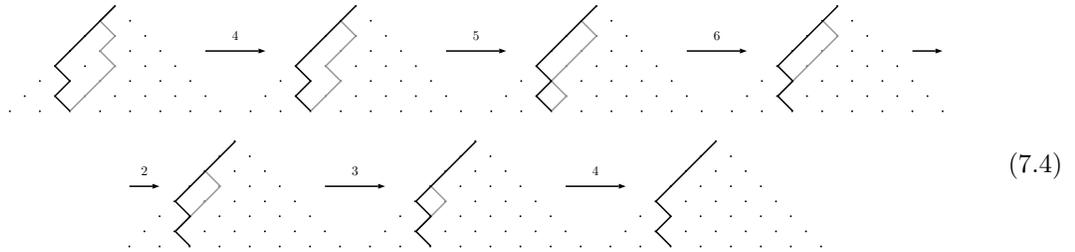
\begin{equation}   \label{eq algorithm}
\DTejem
\end{equation}
Then, $d(\T)=s_4s_3s_2s_6s_5s_4$.
\end{exa}

\begin{defi}
Let $\lambda \in \bip$. $\omega^{\mu_{\lambda}}$ and $\omega_\lambda$ denote the $n$-walks associated with $\mathfrak{t}^{\mu_\lambda}$ and $\T_\lambda$, respectively. The region on the Pascal triangle bounded by $\omega^{\mu_{\lambda}}$ and $\omega_\lambda$ is considered.  This region splits naturally into $q_\lambda +1$ subregions. These subregions are denoted by $B_{\underline{a}}$ ($1\leq a \leq q_\lambda +1$) according to their heights. That is, $B_{\underline{j}} $ is located above $B_{\underline{j+1}}$ for all $j$. Finally, $d(B_{\underline{a}})\in S_n$ is defined as the permutation obtained by performing Algorithm \ref{algorithm} on $B_{\underline{a}}$.
\end{defi}

\begin{exa}
Let $\lambda = (23, 2) \in \operatorname{Bip}_1(25)$, $l = 7$, and $m = 3$.  Under these conditions, $\mu_\lambda = (13,12)$ and  $q_\lambda =2$. The region bounded by $\omega^{\mu_{\lambda}}$ and $\omega_\lambda$ splits into three subregions $B_{\underline{1}},B_{\underline{2}}$, and $B_{\underline{3}}$, as shown in (\ref{lalalalalla}).
Applying Algorithm \ref{algorithm} to these subregions yields
\begin{align*}\label{align subregions}
 d(B_{\underline{1}})& = (s_8s_7s_6s_5)(s_9s_8s_{7})(s_{10}s_9)(s_{11}), \\
d(B_{\underline{2}})&=(s_{15}s_{14})(s_{16}), \\
d(B_{\underline{3}})&= (s_{22}s_{21}s_{20}s_{19})(s_{23}s_{22}s_{21})(s_{24}s_{23}).
\end{align*}
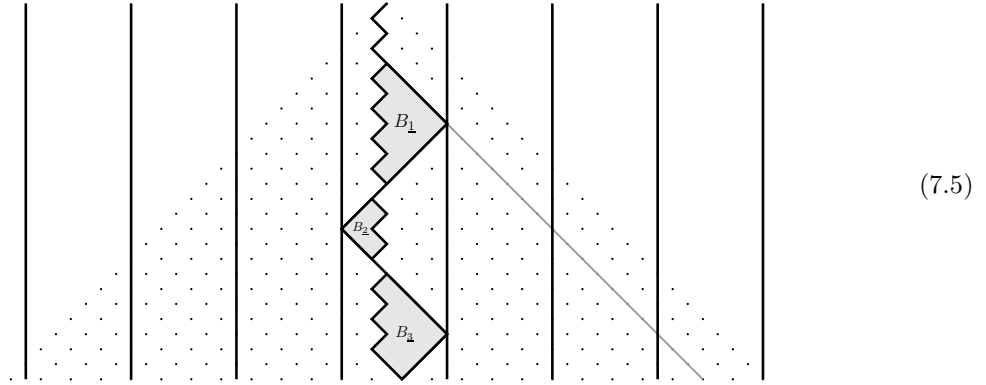
\begin{equation} \label{lalalalalla}
\tribi
\end{equation}
\end{exa}

Step $2$ in Algorithm \ref{algorithm} can be performed in several different ways. Consequently, the obtained reduced expression is not unique.
The following lemma provides a special application of Algorithm \ref{algorithm} to obtain a particular reduced expression for $d(B_{\underline{j}})$. This reduced expression is the key to the proof of Lemma \ref{lema central reflection equal to product of Y's}.

\begin{lem}\label{lembloques}
Let $\lambda\in \bip$.  Then,  $d(\T_{\lambda})\in \mathfrak{S}_n$ can be written as $d(\T_{\lambda})=d(B_{\underline{1}})d(B_{\underline{2}})\cdots d(B_{\underline{q_\lambda+1}})$. Furthermore,  $d(B_{\underline{a}})d(B_{\underline{b}})= d(B_{\underline{b}})d(B_{\underline{a}})$ for all $1\leq a,b \leq q_\lambda +1$. Finally, each $d(B_{\underline{a}})$ can be written as
\begin{equation}\label{forma1}
(s_{k_1}s_{k_1-1}\cdots s_{j_1})(s_{k_2}s_{k_2-1}\cdots s_{j_2})\cdots(s_{k_d}s_{k_d-1}\cdots s_{j_d}),
\end{equation}
for some integers $k_i$ and $j_i$ satisfying  $k_i\geq j_i$, $k_{i+1}=k_i+1$, $j_{i+1}=j_i+2$ and $k_1= \underline{a}$.
\end{lem}

\begin{dem}
It suffices to show that each $d(B_{\underline{a}})$ can be written in the form (\ref{forma1}), but this is immediate by performing Algorithm \ref{algorithm} on the relevant subregion from the bottom to the top, as  shown in (\ref{triangulogeneral}).
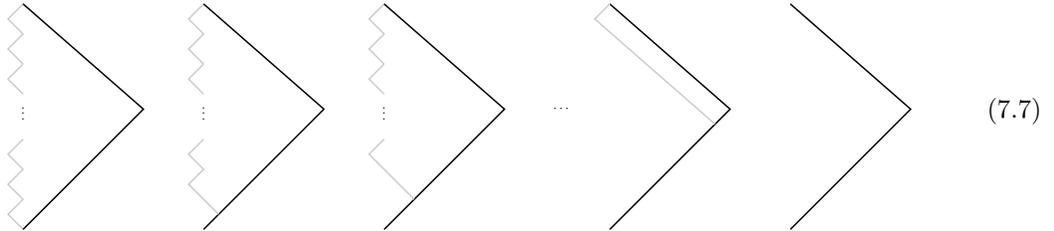
\begin{equation}\label{triangulogeneral}
\demotria
\end{equation}
\end{dem}

Assuming that $d(B_{\underline{a}})$ is one of the reduced expression appearing in the decomposition of $d(\T_{\lambda})$ as in Lemma \ref{lembloques}, we have
\begin{equation}
 k_1-j_1\leq l-3,
\end{equation}
which can be easily deduced from (\ref{condition second equation}). By construction, this implies a more general condition. Namely,
\begin{equation} \label{condicionkj}
 k_a-j_a\leq l-(a+2)
\end{equation}
for all $1\leq a \leq d$. In the following lemma, the notation and the conditions in Lemma \ref{lembloques} are retained.

\begin{lem}\label{lemacom}
Let $1\leq a \leq d$. It is assumed that $j_a\leq i \leq k_a$. Then,
$$2\leq k_a-i+a\leq l-2 $$
unless $a=1$ and $i=k_a$.
\end{lem}

\begin{dem}
This is an immediate consequence of Lemma \ref{lembloques} and (\ref{condicionkj}).
\end{dem}
\par

Lemma \ref{lema central reflection equal to product of Y's} can now be proved. To this end, it is convenient to introduce the following notation. Let $w\in \mathfrak{S}_n$ be a $321$-avoiding permutation. It is assumed that $w=s_{i_1}\cdots s_{i_k}$ is a reduced expression for $w$, and let $\psi_w= \psi_{i_1} \cdots \psi_{i_k} \in b_n $. By (\ref{kl8}), it is clear that $\psi_w$ is well defined. That is, it does not depend on the reduced expression for $w$.

\medskip
\begin{demespecialB}
By Lemma \ref{lembloques}, the subregions $B_{\underline{a}}$ are independent. Therefore, it suffices to prove the result when there exists only one region $B_{\underline{1}}$. Let $w=d(B_{\underline{1}})$. This implies that the following should be proved:
\begin{equation}\label{eq prueba lema}
\psi_{\T_{\lambda}\T_{\lambda}}^{\mu_{\lambda}}=\psi_we(\bi^{\mu_{\lambda}})\psi_{w^{-1}}=\psi_{\underline{1}}^2
e(\bi^{\lambda})=Y_1.
\end{equation}
By using the reduced expression of $w$ given in (\ref{forma1}), we have
 \begin{equation}\label{descom}\psi_{\T_{\lambda}\T_{\lambda}}^{\mu_\lambda}=(\psi_{k_1}\psi_{k_1-1}\cdots \psi_{j_1})\cdots(\psi_{k_d}\psi_{k_d-1}\cdots \psi_{j_d})
e(\bi^{\mu_{\lambda}})(\psi_{j_d} \cdots \psi_{k_d-1}\psi_{k_d})\cdots(\psi_{j_1} \cdots \psi_{k_1-1}\psi_{k_1}).
\end{equation}
Then, the proof is complete if it can be shown that in (\ref{descom}) all $\psi_a$'s on the right-hand side of $e(\bs{i}^{\mu_{\lambda}})$ cancel with their left-hand side counterparts except for the factors $\psi_{k_1}$ located at the ends.  To this end, induction is used on the position of the $\psi_a$'s  in (\ref{descom}). Here, the $\psi_a$'s are arranged from the center to the extremes. For instance, in this order, $\psi_{j_d}$ is the first and $\psi_{k_1}$ is the last.  \par
The first step of the induction is to show that the $\psi_{j_d}$'s located on the right and on the left of  $e(\bs{i}^{\mu_{\lambda}})$ cancel when $j_d \neq k_1$.
By (\ref{kl5}), we have
$$
\psi_{\T_\lambda\T_\lambda}^{\mu_\lambda}=\psi_{ws_{j_d}}\psi_{j_d}e(\bi^{\mu_{\lambda}})\psi_{j_d}\psi_{s_{j_d
}w^{-1}}
=\psi_{ws_{j_d}}\psi_{j_d}^2e(s_{j_d}\bi^{\mu_{\lambda}})\psi_{s_{j_d}w^{-1}}
$$
Let $s_{j_d}\bi^{\mu_{\lambda}}=(r_1,\ldots , r_n)$. Then, by (\ref{kl11}), the $\psi_{j_d}$'s cancel if $|r_{j_d}-r_{j_d+1}|>1$. Indeed, $\mathfrak{S}_n$ acts on $\bi^{\mu_\lambda}$ by permuting the entries. Furthermore, $w\bi^{\mu_{\lambda}}=\bi^{\lambda}$. By considering the reduced expression obtained for $w$ in Lemma \ref{lembloques}, it is clear that $ws_{j_d}(j_d)=k_d+1$ and $ws_{j_d}(j_d+1)=j_d-(d-1)$ are the positions of the residues $r_{j_d}$ and $r_{j_d+1}$ in $\bi^{\lambda}$, respectively. Moreover, if $\bi^{\lambda}=(i_1,\ldots,i_t,\ldots,i_n)$, then $i_{t}-i_{t-1}=1$ for all  $j_1\leq t \leq k_d+1$, because for such a $t$, the $n$-walk associated with $\T^{\lambda}$ is a straight line. Finally, it follows from Lemma \ref{lemacom} that
$$2\leq ws_{j_d}(j_d)-ws_{j_d}(j_d+1)\leq l-2,$$
which implies $|r_{j_d}-r_{j_d+1}|>1$. This proves the first step of the induction.\par

It is now assumed that the first $i-1$ $\psi_{a}$'s in (\ref{descom}) have been cancelled. That is,
$$
\psi_{\T_\lambda\T_\lambda}^{\mu_\lambda}=\psi_{w\sigma_{i-1}}e(\sigma_{i-1}^{-1}\bi^{\mu_{\lambda}})
\psi_{\sigma_{i-1}^{-1}w^{-1}},$$
where $\sigma_{j}$ denotes the product of the first $j$ factors in the reduced expression of $w^{-1}$ obtained by reversing the reduced expression of $w$ in Lemma \ref{lembloques}. Let $\psi_i$ be the $i$-th factor of $\psi_{w^{-1}}$. Then, 
$$\psi_{\T_\lambda\T_\lambda}^{\mu_\lambda}=\psi_{w\sigma_{i}}\psi_{i}e(\sigma_{i-1}^{-1}\bi^{\mu_{\lambda}})
\psi_{i}\psi_{\sigma_{i}^{-1}w^{-1}}=\psi_{w\sigma_{i}}\psi_{i}^2e(\sigma_{i}^{-1}\bi^{\mu_{\lambda}})
\psi_{\sigma_i^{-1}w^{-1}}.$$
 By construction,  $j_a\leq i \leq k_a$ for some $a=1,\ldots,d$. As $w\sigma_{i}(i)=k_a+1$ and $w\sigma_{i}(i+1)=j_a-(a-1)$ are the positions of the residues $r_{i}$ and $r_{i+1}$, respectively, in the tuple $\bi^{\lambda}$, it suffices to repeat the argument of the first step of the induction to finish the proof.
\end{demespecialB}

\section{Vector space isomorphism between $b_n(\lambda) $ and $\al $}

The aim in this section is to prove Conjecture \ref{conjecture one} at the level of $\mathbb{Z}$-graded vector spaces. It is recalled from the introduction that given a choice of the parameters $q$ and $m$, a geometry of alcoves on the real line has been defined. In this context, the generators $\ese $ and $\te $ of $W$ were identified with the reflections in the walls located at $-m$ and $-m+l$, respectively. Moreover, the alcove containing $0$ was denoted by $\mathcal{A}^0$. Furthermore, given $w\in W$, $\mathcal{A}^w$ was defined to be $w \cdot \mathcal{A}^{0}$. Given a positive integer $n$, a bijection is defined as follows:
\begin{equation}
\begin{array}{rcl}
\bip & \longrightarrow & \Lambda_n:=\{-n,-n+2,\ldots , n \} \\
\lambda =(\lambda_1 , \lambda_2) & \longrightarrow & \lambda_1 -\lambda_2 .
\end{array}
\end{equation}
This allows identifying one-line bipartitions with integers. Henceforth, $\lambda \in \bip $ and $w\in W$ are fixed such that $\lambda \in \mathcal{A}^{w}$. The first goal in this section is to find a basis of $b_n(\lambda)$. We recall that $b_n(\lambda) = e(\boldsymbol{i}^\lambda   ) b_n(q,m) e( \boldsymbol{i}^\lambda )$.

\begin{lem}  \label{lema basis bn lambda}
Let $\lambda \in \bip$. Then, the  set
\begin{equation}\label{basis for bn lambda}
  \left\{ \psi_{\mathfrak{st}}^{\mu} \mbox{ }| \mbox{ } \mu \in \bip \mbox{, } \mathfrak{s}, \mathfrak{t} \in \operatorname{Std}(\mu) \mbox{ and }  \bs{i}^\mathfrak{s}=\bs{i}^\mathfrak{t}= \bs{i}^\lambda \right\}
\end{equation}
is a basis of $b_n(\lambda)$.
\end{lem}
\begin{dem}
The result follows by combining (\ref{ortogonality property}), Theorem \ref{teo graded cellular basis blob}, and the definition of $b_n(\lambda)$.
\end{dem}

\par
To prove that  $\al \cong b_n(\lambda)$, as vector spaces, it suffices to construct a bijection between the Double Leaves Basis $\mathbb{DL}_w$ of $\al $ and the basis of $b_n(\lambda)$ described in Lemma \ref{lema basis bn lambda}, which is denoted by $\mathcal{B}_n(\lambda)$.
\par

To control the elements of $\mathcal{B}_n(\lambda)$, all standard bitableaux with the same residue sequence as that of $\mathfrak{t}^\lambda$ should be known. For the present purposes, it is preferable to work with walks on the Pascal triangle rather than bitableaux. Let $\omega^{\lambda}$ denote the $n$-walk associated with $\mathfrak{t}^{\lambda}$. Then, by Lemma  \ref{lemma residues caminos}, the set of all $n$-walks $\omega$ such that $\omega \sim \omega^{\lambda}$ should be determined. Given a walk on the Pascal triangle, a subset of consecutive edges is called a \emph{wall-to-wall
step} if these edges form a straight line between two walls of the same alcove. Then, a walk $\omega $ satisfying $\omega \sim \omega^{\lambda}$ can be described as follows. $\omega $ first coincides with $\omega^{\lambda}$ until the first contact of $\omega^{\lambda}$ with a wall. Subsequently, $\omega $ makes $q_\lambda$ wall-to-wall steps. Finally, $\omega $ is completed with a straight line to the level $n$ in either direction. Figure (\ref{dibujo caminos similares al mayor}) clarifies this description. \par

Thus, to specify a walk $\omega \sim \omega^{\lambda}$, it suffices to describe its behavior on the walls. The focus is on the walk's behavior immediately before/after it touches a wall. Obviously, in both cases, a walk may move toward to or away from the fundamental alcove $\mathcal{A}^{0}$.  Figure (\ref{dibujo rebote murallas}) below shows the four possibilities for a walk when it touches a wall.
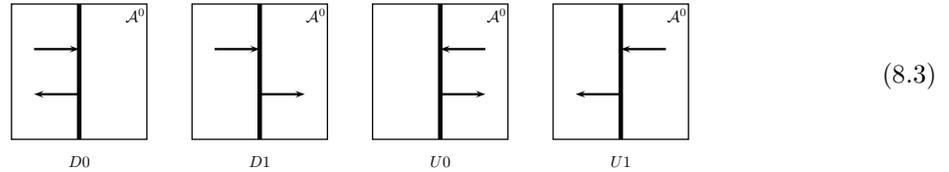
\begin{equation}\label{dibujo rebote murallas}
  \psscalebox{.6}{\rebotesmurallas  }
\end{equation}
 For instance, the leftmost figure in (\ref{dibujo rebote murallas}) shows a walk that, immediately before (resp., after) touching a wall, moves toward to (resp., away from) the fundamental alcove. In this setting, each walk $\omega \sim \omega^{\lambda}$ determines a sequence of elements in $\{D0, D1, U0, U1\}$. It is recalled from the construction of $\mathbb{T}_w$ in Section 2, that a light leaf of $w$ can also be described as a sequence of elements in $\{D0, D1, U0, U1\}$. This allows establishing a bijection between the set of  walks  $\omega \sim \omega^{\lambda}$ and the set of light leaves of $\mathbb{T}_w$.

\begin{lem}  \label{lema bijection walks light leaves}
Let $\lambda \in \bip $ and $w\in W$ such that $\lambda \in \mathcal{A}^{w}$.  Moreover, let $\Lambda_w := \{ x\in W \mbox{ }| \mbox{ }  x\leq w  \}$ and
  $\Lambda_n(\lambda) := \{ \mu \in \bip \mbox{ } |\mbox{ } \exists \mathfrak{t} \in \operatorname{Std}(\mu) \mbox{ such that } \bs{i}^\mathfrak{t}= \bs{i}^\lambda \}$, where $\leq$ denotes the usual Bruhat order of $W$. Given $\mu \in \Lambda_n(\lambda)$, let $\operatorname{Std}_\lambda (\mu) := \{ \mathfrak{t}\in \operatorname{Std}(\mu) \mbox{ } | \mbox{ } \bs{i}^{\mathfrak{t}} = \bs{i}^{\lambda} \}$.  Then, the maps $F$ and $F_\mu$ ($\mu \in \Lambda_n(\lambda )$), defined below, are bijections.
  \begin{enumerate}
    \item $F:\Lambda_n(\lambda) \rightarrow \Lambda_w $, where $F(\mu)=x$ if and only if $\mu \in \mathcal{A}^{x}$ for all $\mu \in \Lambda_n(\lambda)$.
    \item $F_\mu : \operatorname{Std}_\lambda (\mu ) \rightarrow \mathbb{L}_w(F(\mu))$, where for $\mathfrak{t} \in \operatorname{Std}_\lambda (\mu)$, $F_\mu (\mathfrak{t})$ is defined as the unique light leaf in $\mathbb{L}_w(F(\mu))$ such that its sequence of elements in $\{D0, D1, U0, U1\}$ is the same as that of the walk associated with $\mathfrak{t}$.
  \end{enumerate}
\end{lem}
\begin{dem}
This is a straightforward exercise and is left to the reader.
\end{dem}

\begin{cor} \label{coro biyection as vector spaces}
Let $\lambda \in \bip $ and $w\in W$ such that $\lambda \in \mathcal{A}^{w}$. Then, $\al \cong b_n(\lambda)$ as vector spaces.
\end{cor}
\begin{dem}
 It is recalled that a double leaf consists of two light leaves with the same top boundary. Likewise, each element of $\mathcal{B}_n(\lambda)$ is constructed from two standard bitableaux of the same shape and with residue sequence equal to $\bs{i}^{\lambda}$. A map $\mathcal{B}_n(\lambda) \rightarrow \mathbb{DL}_w$ is defined by  $\psi_{\mathfrak{st}}^{\mu} \rightarrow F_\mu (\mathfrak{s})^{a} \circ F_\mu (\mathfrak{t})$ for all $\mu \in \Lambda_n(\lambda)$ and $\mathfrak{s}, \mathfrak{t} \in \operatorname{Std}_\lambda (\mu)$. By using Lemma \ref{lema bijection walks light leaves}, this map is a bijection.
\end{dem}

\begin{rem}
By using the results in \cite[Section 4]{plaza}, it is easy to see that the map defined in the proof of Corollary \ref{coro biyection as vector spaces} is degree preserving. Consequently, $\al \cong b_n(\lambda)$ as $\mathbb{Z}$-graded vector spaces. Moreover, empirical computations suggest that this map is also an algebra isomorphism.
\end{rem}

\end{document}